\documentclass[12pt]{article}
\usepackage{amsfonts,amsmath,latexsym,amssymb}
\usepackage{graphicx}
\usepackage{color}

\newenvironment{proof}{\noindent{{\bf Proof}:\ }}{\vspace{2ex}}

\newenvironment{remark}{{\nin\bf Remark.}}{{\vspace{2ex}}}

\numberwithin{equation}{section}
\allowdisplaybreaks

\newtheorem{theorem}{Theorem}[section]
\newtheorem{lemma}[theorem]{Lemma}
\newtheorem{corollary}[theorem]{Corollary}

\newcommand{\RR}{{\mathbb R}}

\newcommand{\var}{{\mbox{Var}}}

\newcommand{\beas}{\begin{eqnarray*}}
\newcommand{\enas}{\end{eqnarray*}}
\newcommand{\eqs}{\begin{eqnarray*}}
\newcommand{\ens}{\end{eqnarray*}}

\newcommand{\eqa}{\begin{eqnarray}}
\newcommand{\ena}{\end{eqnarray}}
\newcommand{\eq}{\begin{equation}}
\newcommand{\en}{\end{equation}}

\newcommand{\proofbox}{\hspace*{\fill}\mbox{$\halmos$}}
\newcommand{\halmos}{\rule{1ex}{1.4ex}}

\def\ignore#1{}

\def\half{{\textstyle{\frac12}}}

\def\uo{^{(0)}}

\def\bM{{\overline M}}

\def\Ref#1{\eqref{#1}}
\def\a{\alpha}

\def\s{\sigma}
\def\f{\phi}

\def\l{\lambda}
\def\L{\Lambda}

\def\dtv{d_{TV}}
\def\law{{\cal L}}

\def\ep{\hfill $\proofbox$ \bigskip}
\def\Def{\ :=\ }

\def\re{\RR}

\def\giv{\,|\,}

\def\tfrac#1#2{{\textstyle{\frac#1#2}}}

\def\Po{{\rm Po\,}}
\def\BHJ{Barbour, Holst \& Janson}
\def\non{\nonumber}
\def\th{\theta}
\def\bN{{\overline N}}
\def\e{\varepsilon}
\def\m{\mu}

\def\var{{\rm Var\,}}

\def\r{\rho}

\def\t{\tau}
\def\g{\gamma}
\def\h{\eta}
\def\f{\phi}

\def\Bl{\left(}
\def\Br{\right)}
\def\Blm{\left|}
\def\Brm{\right|}

\def\Blb{\left\{}
\def\Brb{\right\}}

\def\ui{^{(1)}}
\def\ut{^{(2)}}

\def\nin{\noindent}
\def\D{\Delta}

\def\ex{{\mathbb E}}
\def\pr{{\mathbb P}}
\def\msk{\medskip}

\def\ff{{\cal F}}
\def\d{\delta}
\def\z{\zeta}
\def\uh{^{(3)}}

\def\Eq{\ =\ }
\def\Le{\ \le\ }

\def\Z{{\mathbb Z}}
\def\C{{\mathbb C}}

\def\bigo{\mathop{{}\mathrm{O}}}

\def\tH{{\widetilde H}}

\def\ignore#1{}

\def\Giv{\,\Big|\,}
\def\cupdot{\cup\kern-8.2pt\cdot\kern5.5pt}

\def\tX{\widetilde X}
\def\tF{\widetilde F}
\def\tG{\widetilde G}

\def\tth{{\tilde\th}}

\def\ggg{\tlaw}

\def\bone{{\bf{1}}}

\def\Def{\ :=\ }

\def\MBP{\bX}

\def\tY{{\widetilde Y}}

\def\hJ{{\widehat J}}
\def\hW{{\widehat W}}

\def\KK{{\cal K}}
\def\htau{{\hat{\tau}}}

\def\p{\pi}

\def\bP{{\overline P}}

\def\bbb{\baL}
\def\mmm{M}

\def\mmc{{\widetilde M}}
\def\nnn{{N}}

\def\MBPK{{\bX^{K}}}
\def\MBPKD{{\bX^{K'}}}

\def\hc{{\hat c}}
\def\nn{{\cal N}}

\def\imag{{\imath}}

\def\tg{{\tilde g}}
\def\rrr{\z}
\def\hX{{\widehat X}}
\def\hA{{\widehat A}}
\def\ha{{\hat a}}

\def\Th{\Theta}
\def\fpp#1{{\ff_{#1}^{++}}}
\def\gda{\g}
\def\tlm{t_{{\rm max}}(\L)}
\def\tml{\tlm}
\def\dbw{d_{{\rm BW}}}
\def\Fbw{F_{{\rm BW}}}
\def\etaL{\h_\L}

\def\vva{{\nu}}
\def\bX{{\overline X}}
\def\hff{{\widehat\ff}}
\def\tff{{\widetilde\ff}}
\def\tlaw{{\widetilde\law}}
\def\bt{{\bar\tau}}
\def\baL{{\overline\law}}
\def\Js{J\uss_s}
\def\ust{^{s,t}}
\def\uss{^{s,s}}
\def\lawtsK{\law^{0,t_s,K}}
\def\bJ{{\overline J}}
\def\bQ{{\overline Q}}

\def\hU{{\widehat U}}
\def\hB{{\widehat B}}
\def\tJ{{\widetilde J}}
\def\cG{{\mathcal G}}
\def\cC{{\mathcal C}}

\begin{document}

\title{A central limit theorem for the gossip process}
\author{
A. D. Barbour\footnote{Institut f\"ur Mathematik, Universit\"at Z\"urich,
Winterthurertrasse 190, CH-8057 Z\"URICH.
Work begun while ADB was Saw Swee Hock Professor of Statistics at the National University of
Singapore, and supported in part by Australian Research Council Grants Nos DP120102728, DP120102398,
DP150101459 and DP150103588.
\msk}
\ and
A. R\"ollin\footnote{Department of Statistics and Applied Probability, National University of Singapore, 
6~Science~Drive~2, 117546 Singapore. Supported in part by NUS Research Grant R-155-000-167-112 and Australian Research Council Grant No.\ DP150101459.
}\\
Universit\"at Z\"urich and National University of Singapore}

\date{}
\maketitle

\begin{abstract}
The Aldous gossip process represents the dissemination of information in geographical space as a process of 
locally deterministic spread, augmented by random long range transmissions.  Starting from a single initially 
informed individual, the proportion of individuals informed follows an almost deterministic path, but for a 
random time shift, caused by the stochastic behaviour in the very early stages of development.  In this paper, 
it is shown that, even with the extra information available after a substantial development time, this broad 
description remains accurate to first order.  However, the precision of the prediction is now much greater, and 
the random time shift is shown to have an approximately normal distribution, with mean and variance that can be
computed from the current state of the process.
\end{abstract}

\makeatletter


\paragraph{MRC subject classification.} 92H30; 60K35, 60J85.

\section{Introduction}\label{intro}

A model for the dissemination of information in space, in which random long-range
contacts facilitate spread, was introduced in Aldous~(2012).  
In an idealized version, proposed by Chatterjee \& Durrett~(2011),
individuals are represented as a continuum, evenly distributed over a two-dimensional torus of 
large area~$L$.
Information spreads locally at constant rate from individuals to their neighbours,
so that a disc of informed individuals, centred on an initial informant, grows steadily in the torus.
However, information is also spread by long range transmissions to other, randomly chosen 
points of the torus, according to a Poisson process, whose
rate is proportional to the area of currently informed individuals.  Any
such transmission  initiates a new disc of informed individuals. 
The process can also be interpreted as a model of the spread of an SI disease,
in which local infection is supplemented by occasional long-range contacts.

With $L_t$ denoting the area of informed individuals by time $t$, Chatterjee \& Durrett~(2011) showed that, after some randomness in
the initial stages of the process, 
the proportion of the torus~$L_t/L$ that has been informed 
by time~$t$ closely follows a particular, deterministic path. The times at which $L_t/L$ increases from almost zero to almost one is relatively short, and occurs around a time $t_L$, which is a fixed multiple of $\log L$. In what follows, we therefore concentrate on times relative to $t_L$. Roughly speaking, Chatterjee \& Durrett~(2011) showed that, for large~$L$, we have
\[
  \frac{L_{t_L+u/\lambda}}{L} \approx \ell(u + U)\qquad \text{for any $u\in\re$}
\]
for some function $\ell$, where~$\lambda$ is a scaling factor related to the speed of spread of information, and where $U$ is a random variable.
The path~$\ell$ is the same
for all realizations of the process, but the position on the path at a particular time varies from
realization to realization because of the random time shift~$U$.
This result was generalized to gossip processes
on rather general homogeneous Riemannian manifolds by Barbour \& Reinert~(2013), hereafter referred to as [BR], as well as 
to related `small world' processes; they also derived a uniform bound on the approximation
error.  In addition, the 
equation describing the deterministic development was interpreted in terms of the
Laplace transform of the limiting random variable corresponding to an associated 
Crump--Mode--Jagers (CMJ) branching process (Jagers, 1975).

By analogy with the theory of Markov population processes (Kurtz 1970, 1971), one might expect that the
fluctuations around the deterministic path of the proportions informed would be approximately
Gaussian, with standard deviation~$O(L^{-1/2})$, at least while the proportion informed is not
too small or too close to~$1$.  Here, however, the random quantity of most interest ---
the difference between the actual course of the process and a prediction of the course based on
information available early in its development --- involves the fluctuations of the process while the 
proportion informed is rather small, and the standard analogy does not apply. 
Instead, in view of the approximation already established, it seems reasonable at times $v \ll t_L$ to predict
the value of $L_{t_L + u/\lambda}/L$ by $\ell(u + \hU(v))$, where~$\hU(v)$ is the expected value of~$U$,
given the information at time~$v$, and to augment the point prediction with a confidence interval
around~$\ell(u + \hU(v))$,
derived from the (approximate) conditional distribution of $L_{t_L + u/\l}/L$, given the current information.

The validity of the procedure is justified in detail in Section~\ref{CLT}.  The broad argument is to
exploit the fact that $L_{t_L + u/\l}/L$ is the probability that a point~$K$, chosen independently and
uniformly at random in~$\cC$, belongs to the informed set~$\law_{t_L + u/\l}$:
\[
     L_{t_L + u/\l}/L \Eq \pr[K \in \law_{t_L + u/\l} \giv \law_{t_L + u/\l}].
\]
As it stands, this changes nothing.  However, it indicates that a good approximation might be
obtained by replacing $\pr[K \in \law_{t_L + u/\l} \giv \law_{t_L + u/\l}]$ by $\pr[K \in \law_{t_L + u/\l} \giv \law_{s}]$,
or, equivalently, replacing $L_{t_L + u/\l}/L$ by~$\ex\{L_{t_L + u/\l}/L \giv \law_s\}$,
for $s < t_L + u/\l$ chosen so that~$s$ is close enough to~$t_L+u/\l$.  In particular, for prediction
from~$v$,  we need to choose~$s \in (v,t_L+u/\l)$ so that
\eq\label{approx-1a}
   \ex_v \bigl|(L_{t_L + u/\l}/L) - \ex\{L_{t_L + u/\l}/L \giv \law_s\}\bigr| \ \ll\ {\rm SD}_v(L_{t_L + u/\l}/L),
\en
where $\ex_v$ and ${\rm SD}_v$ denote expectation and standard
deviation given the information at time~$v$.

The advantage of using $\ex\{L_{t_L + u/\l}/L \giv \law_s\}$ is that $\pr[K \in \law_{t_L + u/\l} \giv \law_{s}]$
can be approximated as the probability of at least one of many small balls, with centres chosen independently
and at random in~$\cC$,
intersecting~$\law_s$.  These balls are the islands in an independent `backwards' gossip process, run for a 
length of time $t_L+(u/\l) - s$ from~$K$.  There are many such balls if~$t_L+(u/\l) - s$ is not too small, and the
intersection probability can be approximated by a Poisson probability, using the Stein--Chen method;
see Lemma~\ref{AB-first-representation}.
The mean of the Poisson distribution can, with considerable effort, be shown to be
close to $\ell(\log[C W(s,v)] + u)$, where~$W(s,v)$ is a quantity that can be simply expressed in terms
of a carefully chosen branching process, and~$C$ is a constant.  Now, given the information available at 
time~$v$, the quantity~$W(v,v)$
(which loosely corresponds to $\exp\{\hU(v)\}$)
is known, and the conditional distribution of the difference $W(s,v) - W(v,v)$
is approximately normal, as is shown in Theorem~\ref{AB-branching-appxn}
in Section~\ref{branching}.  This, in turn, leads to a normal approximation for the difference between
$\ell(\log[C W(s,v)] + u)$ and its prediction $\ell(\log[C W(v,v)] + u)$ at time~$v$. This implies the
main result of the paper, that
\eq\label{rough-theorem}
   \sigma^{-1}\bigl(L_{t_L + u/\l}/L - \ell(\log[C W(v,v)] + u)\bigr) \ \approx_d\ \nn(0,1),
\en
for suitable choice of the standard deviation~$\s$ depending on $u$ and $W(v,v)$; a precise statement is given in Theorem~\ref{ADB-main-theorem}.  
The error in the normal approximation is shown to be small if the {\it number\/} of 
individuals informed at time~$v$ is large, even if their proportion in the whole population may be very small.
For practical purposes, in an epidemic, the very earliest development
may well pass almost unnoticed --- the origins are often obscure --- but prediction on the basis of the information
gained from the first few hundred cases is an important public health goal, in which case using the normal
approximation is reasonable.

\subsection{Detailed formulation}
We now describe the problem in more detail.
We consider the gossip process~$(\law_t,\,t\ge0)$ evolving on a smooth closed homogeneous 
Riemannian manifold~$\cC$ of dimension~$d$,
such as a sphere or a torus, having large finite volume~$|\cC| =: L$ with respect to its intrinsic
metric. An individual at point $P\in \cC$ informed at time~$0$ gives rise to deterministic local spread that informs
the set~$\KK(P,s)$ by time~$s > 0$;  in addition, random `long range transmissions' 
to independent and uniformly distributed points of~$\cC$ occur at rate~$\r$ times the intrinsic volume 
of the set currently informed.  Thus the process
can be constructed from knowledge of the points $0 = \t_0 < \t_1 < \cdots$
of a point process~$\Pi$ on~$\re_+$ (characterized immediately below), together with an independent sequence of independent points 
$P_1,P_2,\ldots$, uniformly distributed in~$\cC$, and an initial point~$P = P_0$. The informed set and its volume are denoted by
\eq\label{ADB-informed-set}
     \law_t \Def \bigcup_{j\colon \t_j \le t} \KK(P_j,t-\t_j)
     \quad\text{and}\quad L_t \Def |\law_t|.
\en
The point
process~$\Pi$ is simple, and has conditional intensity $\r L_t$ at time~$t$ with respect to the
filtration~$(\ff_t,\,t\ge0)$, where $\ff_t := \s((\t_j,P_j),\,j\ge0,\t_j \le t)$.

The sets~$\KK(P,s)$ are assumed to be closed balls, centred at~$P$ and of radius~$s$,  with respect to
a metric that makes~$\cC$ a geodesic space: $P' \in \KK(P,2t)$ exactly when 
$\KK(P,t) \cap \KK(P',t) \neq \emptyset$. 
Since $\cC$ is assumed to be homogeneous, the volume of $\KK(P,s)$ is independent of $P$, and we will therefore denote 
it by $\vva_s = \vva_s(\KK)$.
The sets~$\KK(P,s)$ are also assumed to be locally almost Euclidean in the sense that $\vva_s\approx s^d\vva $ for 
some constant $\vva = \vva(\KK) > 0$. More precisely, we will assume that, for constants $c_g,\g_g>0$, 
\eq\label{AB-volume-approx}
  \bigg|\frac{\vva_s}{s^d\,\vva } - 1\bigg| \Le c_g\bigg(\frac{s^d\,\vva }{L}\bigg)^{\g_g/d},\quad s > 0.
\en 
The quantity~$\vva  > 0$ has physical dimensions $(\mbox{length}/\mbox{time})^d$,
so that $\vva ^{1/d}$ can be interpreted as a local velocity of spread of information in any particular direction.
Assumption~\Ref{AB-volume-approx} is satisfied, for instance, for balls with respect to geodesic 
distance on 
the surface of a $(d+1)$-dimensional sphere of large radius~$R$, when $L = c_d R^d$ and
\[
     \frac{\vva_s}{s^d\,\vva } - 1 \Eq \frac{dR^d}{s^d}\,\int_0^{s/R} (\sin t)^{d-1}\,dt - 1
                     \Eq O\bigl( (s/R)^2 \bigr),
\]
(Li, 2011), in which case we can take
$\g_g = 2$ in all dimensions $d \ge 2$.

Using~\Ref{AB-volume-approx}, the probability of there being no long range transmission before time~$u$ is given by
\[
   \exp\Bigl\{ - \int_0^u \r\vva_s\, ds\Bigr\} \ \approx\ \exp\Bigl\{ - \int_0^u \r s^d\vva \,ds\Bigr\}
                    \Eq \exp\bigl\{ -\r\vva  u^{d+1}/(d+1)\bigr\},
\]
so that the mean time to the first long range transmission is approximately
\[
    \int_0^\infty \exp\bigl\{-\r \vva  u^{d+1}/(d+1)\bigr\} \,du 
             \Eq (\r\vva)^{-1/(d+1)} \int_0^\infty e^{-w^{d+1}/(d+1)}\,dw.
\]
Thus
\eq\label{averyimportantdefinition}
  \l \Def (\r d!\vva )^{1/(d+1)},
\en
having physical dimensions $(1/\mbox{time})$, 
is such that $1/\l$ represents the time scale for the first long range transmission, and then
$\l^{-d}\vva $ reflects the size of the initial neighbourhood when the first long range
transmission occurs; the exact specification of~$\l$ is to make it equal to the growth
rate of the associated CMJ process ([BR], p.986).  For our approximations to be good, the
size of the initial neighbourhood when the first long range transmission occurs should be small
compared to~$L$, so that, defining
\eq\label{AB-Lambda-def}
  \L \Def L \l^d / \vva ,
\en
a quantity without physical dimension, we shall take~$\L$ to be large.  Note that, if this is so,
the approximations made above have small error, in view of~\Ref{AB-volume-approx}.

To start with, the points of~$\Pi$
closely match the birth events of a CMJ process~$\bX$, whose birth intensity as a function
of age~$s$ is given by $\r \vva_s$.  
In fact, the approximation~$\baL_t$ of~$\law_t$, constructed by using the CMJ process~$\bX$ to approximate~$\Pi$ and
with the same sequence of points $(P_j,\,j\ge1)$, is excellent
for times $t \le \a \l^{-1}\log\L$ if $\a < 1/2$ ([BR], \S2.2), and still gives an approximation to the
volume $L_t$ of~$\law_t$ at time~$t$ that is accurate to the first order if $\a < 1$ 
([BR], Theorem~3.2 and (2.23)).
This CMJ approximation takes the form
\eq\label{CMJ-approx}
   {L_t}/L \ \sim\ K\L^{-1}e^{\l t + \log W},\quad t\to\infty,
\en  
for a constant~$K$, where~$W$ is a limiting random variable associated with the CMJ process~$\bX$.
Taking 
\eq\label{defoftsublambda}
  t \Eq t_\L(u) \Def \l^{-1}(\log\L + u),
\en
with~$u \le (\a-1)\log \L$  {\it large and negative\/} in
the range in which this approximation holds, this implies that $L_{t_\L(u - \log W)}/L$
closely follows the curve
\eq\label{ell_0}
    u \ \mapsto\ \ell_0(u),
\en
where $\ell_0(u) := K e^u$.
  
In [BR], Theorem~3.2, an analogous approximation 
\[
     L_{t_\L(u - \log W)}/L \ \approx\ \ell(u + \log\hc_d)
\]
is established, with uniformly small error, for {\it all\/} values of~$u$,
with~$\hc_d$ defined
before~\Ref{AB-min-r-eta-def},
and with the time shift~$U$ given by $\l^{-1}\log W + c$, for a suitably chosen constant~$c$. Clearly, to be compatible 
with~\Ref{ell_0}, $\ell(u) \sim K e^u$ as $u\to -\infty$, as follows from ([BR], following~(2.23)).

For any fixed~$u$, the distribution of $L_{t_\L(u)}/L$ is close to
that of $\ell(u + \log W + \log\hc_d)$, and is a bounded random variable.
Hence it can only be approximately normally distributed, after appropriate centring and
normalization, in circumstances in which the distribution of~$\log W$ is concentrated close to
some fixed value.  This is not true of the distribution of~$W$ at time~$0$. However, when predicting 
from a time $v = \a\l^{-1}\log\L$ for any fixed~$\a$, $0 < \a < 1$, the conditional distribution
of~$W$, given the information up to time~$v$, is concentrated close to
an approximation $W(v,v)$ provided only that $\a > 0$, even though the size of the informed set is still 
relatively small
when compared to~$L$ for any $\a < 1$.  The aim is now to show that
the difference $\D(v) := W(v,v) - W$, suitably normalized, is approximately normally distributed.

It turns out to be easier to work with a `flattened' CMJ process~$\hX$, 
rather than with the original CMJ process~$\bX$.  The process~$\hX$ has birth rate at age~$s$ 
given by $\r s^d\vva $, and is thus the same process for all~$L$, whereas~$\bX$ depends implicitly on~$L$
through the function~$\vva_s$. The quantity~$\l$ then turns out to be the Malthusian 
parameter of~$\hX$.  In a CMJ process with Malthusian parameter~$\m$, 
at large times, a randomly sampled individual has average age approximately~$1/\m$.  For~$\hX$,
$\m = \l$, and
replacing $s$ by~$1/\l$ in~\Ref{AB-volume-approx} confirms that the two CMJ processes $\bX$ and~$\hX$ have
birth rates that are close to each other if~$\L$ is large.  The essentials of the
proof of the normal approximation to~$\D(v)$ are carried out in Section~\ref{branching}.
The argument hinges on examining a collection of (complex valued) martingales~$(W_j(\cdot),\,0 \le j\le d)$ 
associated with~$\hX$, that are defined in~\Ref{AB-W-def} below.
In particular, $W(t,v) := W_0(t)$, $t\ge v$, is non-negative and 
square integrable, having limit $W_0(\infty) =: W$.  It is then shown that $W_0(v) - W$,
suitably normalized, is close enough
to the integral of a function $f(W_0(v),u)$ with respect to an independent standard Brownian motion~$B(u)$,
giving the normal approximation.

The arguments in Section~\ref{CLT}, as outlined before~\Ref{rough-theorem},
rely heavily on comparisons between birth and growth
processes.  The actual process $(\law_t,\,t\ge0)$ is compared with the branching approximation~$\bX$, and
$\bX$ is compared to its flattened version~$\hX$.
Further (flattened) CMJ processes $\hX^+$ and~$\hX^-$ are then introduced, to act as upper and lower bounds for~$\bX$;
the comparison is formalized in Lemma~\ref{AB-dominating-processes}.  All the detailed computations in Section~\ref{CLT} are
made using these processes, including the reduction of the intersection probability in Lemma~\ref{AB-first-representation}
to a tractable form in Lemma~\ref{AB-W-to-Wsv}.

To state our theorem, we take
\eq\label{ADB-W-star-def-new}
   \hW(v) \Def e^{-\l v}\sum_{l=0}^d \sum_{j\in \hJ_v} \frac{\{\l(v-\t_j)\}^l}{l!}
\en
as an approximation to~$W$,
where the set~$\hJ_v$ indexes the set of all {\it non-intersecting\/} neighbourhoods of~$\law_v$.  For each of these,
the radii $(v - \t_j)$ can be determined, and so~$\hW(v)$ can be derived from~$\law_v$.  Then
let $\hc_d := d!/(d+1)$, and
\eq\label{AB-min-r-eta-def}
   \rrr(d) \Def 
       \begin{cases}
         1/2  &\text{if $d \le 6$,}\\
         1 - \cos(2\pi/d)  &\text{if $d \ge 7$,}
        \end{cases}
\en 
and define 
\eq\label{ADB-l-def}
   \ell(u) \Def 1 - \f_\infty(e^u), \quad\mbox{where}\quad \f_\infty(\th) \Def \ex\{e^{-\th W}\},
\en
where~$W$ is as above; see also~\Ref{AB-W-def} and~\Ref{W-convergence}.
Let~$\dbw$ denote the bounded Wasserstein distance between probability measures on~$\re$:
\[
   \dbw(P,Q) \Def \sup_{f\in \Fbw}\Bigl\{\Bigl|\int f\,dP - \int f\,dQ \Bigr|\Bigr\},
\]
where~$\Fbw$ consists of all Lipschitz functions~$f\colon \re\to [-1,1]$ whose Lipschitz constant
is at most~$1$. 
The theorem is as follows.

\begin{theorem}\label{ADB-main-theorem}
With the above definitions,
suppose that $v = \a\l^{-1}\log\L$ for $0 < \a < 2\min\{\g_g/d,\rrr(d)/(1 + \rrr(d))\}$,
where~$\g_g$ is as in~\Ref{AB-volume-approx}. 
Then, for any $u_1 < u_0 \in \re$, there exists a $\g > 0$ 
and an event $E^*(v) \in \s(\law_v)$  with $\pr[E^*(v)^c] = \bigo(\L^{-\g})$ such that
\eqs
    \lefteqn{
       \dbw\bigl(\law\bigl\{e^{\l v/2}\{L_{t_\L(u)}/L - \ell(u + \log[\hc_d \hW(v)])\} \,\big|\,\ff_v \cap E^*(v)\bigr\},
             \nn(0,\s^2(u,\hW(v))) \bigr) } \\ 
        &&\Eq \bigo(\L^{-\g}), \phantom{XXXXXXXXXXXXXXXXXXXXXXXXXXXXXXXXX}
\ens
uniformly in $u_1 \le u \le u_0$, where  $t_\L(u) = \l^{-1}(\log\L + u)$ as in \eqref{defoftsublambda} and
\[
     \s^2(u,w) \Def \frac{\{D\ell(u + \log[\hc_d w])\}^2}{(d+1)w}.
\]
\end{theorem}

\nin 
So, for instance, for spherical neighbourhoods in $d \le 6$, it is possible to take any~$\a$
strictly between $0$ and~$2/3$ in Theorem~\ref{ADB-main-theorem}.
The order statements can be replaced by inequalities, valid for all~$\L$ sufficiently large, in which the 
constants depend only on~$d, u_1$ and~$u_0$;  however, the
lower bound on the value of~$\L$ then also involves~$\a$ and the constants $c_g$ and~$\g_g$ from \Ref{AB-volume-approx}.

In fact, the proof shows a little more: that we could realize the normal random
variables $\nn(0,\s^2(u,W(v,v)))$, for different values of~$u$, as $\s(u,W(v,v))N$ for the {\it same\/} standard normal
random variable~$N$.
The interpretation of this is that the fluctuations in $L_{t_\L(u)}/L$ are essentially those of 
$\ell(u + \log[\hc_d W])$, and that the remaining randomness after time~$v$ is overwhelmingly that of
the difference $W - W(v,v)$, a single random variable.  
This, at first sight surprising, result reflects the phenomenon common to branching processes, that the 
randomness determining the growth of a super-critical branching process occurs at the very 
beginning of its development.

\section{The branching process}\label{branching}

In this section, we investigate the limit~$W$, as $t\to\infty$, of a martingale~$W(t)$ associated with a 
particular CMJ branching process.  We show that $(W(t)-W)$ is approximately 
normally distributed, and give an explicit bound on the accuracy of the approximation.  Although,
for a (multitype) Galton--Watson process, a central limit theorem of this sort is not
difficult to establish (Asmussen \& Hering, Theorem~7.1), the corresponding theorems for 
general CMJ processes seem not
to be available. Here, we are able to exploit the particular structure of our CMJ process to prove what
we need.

We start by identifying the branching process that we work with, which can be expressed as a Markov
process in a $(d+1)$-dimensional space.  The properties of the coordinate processes $(H_j(t),\,0\le j\le d)$,
and of some equivalent (complex valued) martingales $(W_j(t),\,0\le j\le d)$ are established in 
Lemma~\ref{ADB-martingales}.  The component~$W_0$ is a non-negative real valued martingale, and~$W$ is
its limit as $t\to\infty$.  Using Kolmogorov's inequality, the fluctuations of the sample paths of the 
processes~$W_j$ are controlled in Lemma~\ref{AB-W-bounds}, and this in turn gives control over the
processes $H_j$.

The martingale difference $W_0(v+t) - W_0(v)$ is written in~\Ref{AB-W-difference} as an integral of an 
explicit function of the process
$H_{d+1}(u) := \l\int_0^h H_d(w)\,dw$ with respect to a standard compensated Poisson process.
Using the control that we have over the~$H_j$, we determine successively simpler approximations
to this process, in \Ref{AB-X1-def} and~\Ref{AB-X2-def}, at each stage making sure that the error
incurred is sufficiently small (Lemma~\ref{AB-smooth-integrand} and Corollary~\ref{AB-X1X2-cor}).  
Finally, in~\Ref{AB-X3-def}, 
an expression is obtained in which integration with respect to the 
compensated Poisson process has been replaced by integration with respect to standard Brownian motion,
and this can be used with an error controlled in Lemma~\ref{AB-KMT}.  The results of these
steps are collected as a functional approximation in Theorem~\ref{AB-branching-appxn}. 
The version that is used to prove Theorem~\ref{AB-main-theorem} in Section~\ref{CLT} is given
as Corollary~\ref{AB-CLT-cor}.

\subsection{Properties of the flattened process}
The first step is to determine a suitable~$W$.  We do so
by way of a `flattened' version~$\hX$ of the CMJ branching process~$\bX$.  The process~$\hX$ 
is the counting process associated with a
point process $(\htau_j,\,j\ge0)$ on~$\re_+$, with $\htau_0=0$ a.s., whose compensator is
given by $\hA(t) := \int_0^t \ha(u)\,du$, where $\ha(u) := \r \vva  \sum_{j:\,\htau_j \le u} (u-\htau_j)^d$,
and where~$\r$, as before, denotes the intensity per unit volume. 
At time~$t$, $\hX(t)$ can be thought of as consisting of $M_0(t) := 1 + \max\{r\geq 0\colon\,\htau_r \le t\}$ 
neighbourhoods, whose volumes at time~$t$ are given by~$(t - \htau_r)^d \vva $,  asymptotically close to,
but not the same as the volume~$\vva_{t - \htau_r}$.  The intensity~$\ha$ is then precisely that
of a CMJ process, in which neighbourhoods play the part of individuals, and the point process~$\xi$ of an
individual's offspring is an inhomogeneous Poisson process with rate $\r\vva s^d$ at age~$s$. The
mean number of offspring of an individual is thus infinite, but the Malthusian parameter~$\l$,
chosen so that the equation
\[
    \int_0^\infty e^{-\l s} \r \vva  s^d\,ds \Eq 1 
\]
is satisfied, is finite, and is given by $\l := (d!\r \vva )^{1/(d+1)}$.  Note that
\eq\label{xi-standard-process}
       (\xi(t),\,t\ge0) \ =_d\ (\xi^1(\l t),\,t\ge0),
\en
where~$\xi^1$ is the inhomogeneous Poisson process with rate $s^d/d!$ at age~$s$.

We can immediately deduce some useful general properties of the process~$\hX$.  To start with, because
the variance of the discounted offspring number $\int_0^\infty e^{-\l s}\xi(ds)$ is finite,
being given by $\int_0^\infty e^{-2\l s} \r \vva  s^d\,ds$, it follows 
from Ganuza \& Durham~(1974, Theorem~1) that there exist finite constants $c_1$ and~$c_2$ such that, for all $u>0$,
\eq\label{AB-N-star-moments}
      e^{-\l u}\ex M_0(u) \Le c_1;\qquad e^{-2\l u} \ex\{ M_0^2(u)\} \Le c_2;
\en
in view of~\Ref{xi-standard-process}, $c_1$ and~$c_2$ depend only on~$d$.
Then the intensity~$\ha(u)$ can be expressed as $\r \vva  M_d(u)$, where
\eq\label{AB-L-star-def}
    M_d(u) \Eq  u^d + \int_{(0,u]}(u-v)^d M_0(dv) \Eq d\int_0^u  (u-v)^{d-1}M_0(v)\,dv.
\en
This in turn implies from~\Ref{AB-N-star-moments} that
\eq\label{AB-L-star-moments}
   e^{-\l u} \ex M_d(u) \Le c_1 d!\l^{-d};\qquad e^{-2\l u} \ex\{M_d^2(u)\} \Le c_2\{d!\l^{-d}\}^2,\qquad u > 0,
\en
using Cauchy--Schwarz for the second inequality.

However, $\hX$ also has special structure that will prove useful in what follows,
relating to the sums
\eq\label{AB-M-defs}
   M_l(t) \Eq \sum_{j=1}^{M_0(t)} (t-\t_{j-1})^l,\quad l\ge1,
\en
of the $l$-th powers of the ages of the neighbourhoods. Note that~$M_d(t)$ is as defined previously, and that
\eq\label{AB-M-eqns}
   \begin{array}{rl}
   \dfrac{d}{dt} M_1(t) &=\ M_0(t) \quad \mbox{for a.e. }t;  \\[2ex]
   \dfrac{d}{dt} M_i(t) &=\ i M_{i-1}(t), \quad i\ge 2.
   \end{array}
\en
Since $M_0$ has intensity $\ha = \r \vva  M_d$, letting~$Z$ denote a unit rate Poisson process,
we can write
\eq\label{AB-incidence-vol}
   M_0(t) \Eq  M_0(0) + Z\left(\r \vva  \int_0^t M_d(u)\,du \right).
\en
Defining $H_i(t) := M_i(t)\l^i/i!$, for any $\l > 0$,
the equations~\Ref{AB-M-eqns} reduce to
\eq\label{AB-H-eqns}
  \begin{array}{rl}
   \dfrac{d}{dt} H_1(t) &=\ \l H_0(t) \quad \mbox{for a.e. }t; \\[2ex]
   \dfrac{d}{dt} H_i(t) &=\ \l H_{i-1}(t), \quad i\ge 2;
  \end{array}
\en
with the particular choice $\l := (d!\r \vva )^{1/(d+1)}$, equation~\Ref{AB-incidence-vol} becomes
\eq\label{AB-kappa-vol}
    H_0(t) \Eq M_0(0) + Z\left(\r \vva  \int_0^t d! \l^{-d} H_d(u)\,du\right) \Eq  H_0(0) + Z(H_{d+1}(t)),
\en
so that $\hA(t) = H_{d+1}(t)$. In particular, from \Ref{AB-H-eqns} and~\Ref{AB-kappa-vol},
it follows that the process~$\tH$ defined by
\eq\label{ADB-tH-is-Markov}
   \tH(t) \Def (H_0(t),H_1(t),\ldots,H_{d}(t))
\en
is a Markov process. It also follows directly from \Ref{AB-H-eqns} and~\Ref{AB-kappa-vol}, or as
a consequence of~\Ref{xi-standard-process}, that
\eq\label{AB-standard-H}
  \{\tH(t),\,t\ge0\} \ =_d\ \{\tH^1(\l t),\,t\ge0\},
\en
where $\tH^1$ denotes the process with $\l = 1$.  Note that~$\r$ may depend on~$L$, as also may~$\l$.

In order to describe the properties of the process~$\hX$ in more detail, we introduce the (complex valued)
processes
\eq\label{AB-Wj-def}
   W_j(t) 
               \Eq 1 + \int_{(0,t]} e^{-\l x_ju}\{ M_0(du) - \hA(du)\},
\en
where $x_j := \exp\{2\pi \imag j/(d+1)\} \in \C$,  $j\in \{0,1,\ldots, d\}$, which are martingales
with respect to the natural filtration $(\hff_t,\,t\ge0)$ of~$\hX$.
In particular, for $j=0$, we have $x_j=1$, and
\eq\label{AB-W-def}
    W(t) \Def W_0(t) 
             \Eq 1 + \int_{(0,t]} e^{-\l u}\{ M_0(du) - \hA(du)\}
\en
is a real valued, c\`adl\`ag martingale, and plays a key part our arguments.  It is shown in the
next lemma that it is also non-negative, and the rest of the section is then devoted to proving
a normal approximation to $e^{\l t/2}(W(t) - W(\infty))$, which is the basis for the central limit
theorem for the gossip process itself.
Note that the distribution of~$W(\cdot)$ can be derived from the corresponding martingale~$W^1(\cdot)$
for the process with $\l = 1$, since, from~\Ref{AB-standard-H}, 
\eq\label{AB-standard-W}
        \{W(t),\,t\ge0\}\ =_d\ \{W^1(\l t),\,t\ge0\};
\en
from this, it also follows that the distribution of~$W(\infty)$ is the same for all~$\l$.
The remaining martingales~$W_j$ are useful, because they enable the quantities~$H_j(\cdot)$ to be
expressed in a tractable form, as in the next lemma.

\begin{lemma}\label{ADB-martingales}
With notation as above, we have
\[
     W(t) \Eq \sum_{r=0}^d  e^{-\l  t}H_r(t) \ >\ 0,\quad t \ge 0,
\]
and
\[
   e^{-\l t}H_j(t) \Eq \frac1{d+1} \sum_{l=0}^d x_j^l\, e^{-\l(1- x_l) t}\, W_l(t).
\]
\end{lemma}

\proof
It follows from~\Ref{AB-H-eqns} that, for any~$x \in \C$, 
\eqs
    \frac d{dt} \{ e^{-\l xt}x^r H_r(t) \} &=& \l x e^{-\l xt}\{ - x^r H_r(t) + x^{r-1}H_{r-1}(t) \},  
   \quad r\ge1,
\ens
and, by partial integration, that
\[
     \int_{[0,t]} e^{-\l xu} H_{0}(du) \Eq e^{-\l xt} H_{0}(t) + \l x \int_0^t e^{-\l xu} H_{0}(u)\,du.
\]
Hence
\[
   \frac d{dt} \sum_{r=1}^d \{ e^{-\l xt}x^r H_r(t) \} \Eq \l x e^{-\l xt}\{ - x^d H_d(t) + H_0(t) \},
\]
and thus
\eq\label{AB-Wx-def}
    \sum_{r=0}^d \{ e^{-\l xt}x^r H_r(t) \} 
             \Eq \int_{[0,t]} e^{-\l xu}\{H_0(du) - \l x^{d+1} H_d(u)\,du\}.
\en
Taking $x=x_j$ for any~$j \in \{0,1,\ldots,d\}$, we have $x^{d+1} = 1$, 
making the right hand side equal to~$W_j(t)$, because $\l H_d(u)\,du = H_{d+1}(du) = \hA(du)$,
by \Ref{AB-H-eqns} and~\Ref{AB-kappa-vol}; hence
\eq\label{ADB-Wj-repn}
   W_j(t) \Eq \sum_{r=0}^d \{ e^{-\l x_j t}x_j^r H_r(t) \}.
\en
The first statement of the lemma follows by taking $j=0$, and the second by using the orthogonality 
relation $\sum_{l=0}^d x_j^l x_l^r = (d+1)\d_{jr}$.
\ep

Now, writing $r_{j} := \Re{x_j}$  and noting that $\ha(u) = \l H_d(u) \le \l e^{\l u}W(u)$,
it follows from~\Ref{AB-Wj-def} that, for $0\le j\le d$ and for $v < t < w$,
\eqa\label{AB-MG-variance}
  \ex\{|W_j(w) - W_j(t)|^2 \giv \hff_v\} &=& \int_{(t,w]} e^{-2\l r_{j} u}\ex\{\ha(u) \giv \hff_v\}\,du \non\\
      &\le& W(v) \int_{(t,w]} \l e^{-\l (2r_{j} - 1) u}\,du.
\ena
Using this bound with $v=0$, we see that the variances of the terms with $1 \le l \le d$
in the sum in Lemma~\ref{ADB-martingales} converge to zero as $t \to \infty$.  However, the
term with $l=0$ remains significant as $t\to\infty$, since, by~\Ref{AB-MG-variance} with $v=0$ and $j=0$,
it follows that $W(\cdot)$ is square integrable, and that 
\eq\label{W-convergence}
      W \Def W(\infty) \Def \lim_{t\to\infty} W(t)\quad \mbox{exists a.s.};\quad  \mbox{and }  \ex W \Eq 1,\quad \var W \le 1. 
\en 
Note that the distribution of~$W$, through its Laplace transform~$\f_\infty$ as in~\Ref{ADB-l-def}, already
appears in the statement of Theorem~\ref{ADB-main-theorem}, and is the same for all~$\l$,
as remarked following~\Ref{AB-standard-W}.
Thus each of the~$H_j$ satisfies 
\eq\label{ADB-H-asymp}
      e^{-\l t}H_j(t)\ \to_P\ W/(d+1)\quad\mbox{as}\ t\to\infty.  
\en
We shall exploit more detailed versions of these asymptotics in Section~\ref{CLT}.

In order to use Lemma~\ref{ADB-martingales} to describe further the behaviour of the~$H_j(t)$, 
we need good control of the fluctuations of the processes~$(W_l,\, 0\le l\le d)$.
As indicated by~\Ref{AB-MG-variance}, their asymptotic behaviour depends substantially
on whether or not $r_l > 1/2$.  Note, for future reference, that 
$\min\{(1-r_1),1/2\} = \rrr(d)$, where $\rrr(d)$ is as in~\Ref{AB-min-r-eta-def}.

\begin{lemma}\label{AB-W-bounds}
 For any $1\le l\le d$ and $0 < \h < \min\{(1-r_l),1/2\}$, and for any $K>0$, define the events
\eqs
   E_{1l}^\h(v;K) &:=& \Bigl\{ \sup_{t\ge v} \{e^{-\l t(1-r_l-\h)}|W_l(t) - W_l(v)|\} \le K \Bigr\};
\ens
similarly, for $0 < \h < 1/2$, define
\[
   E_{10}^\h(v;K) \Def \Bigl\{\sup_{t\ge v} \{e^{\l\h t}|W(t) - W(\infty)|\} \le K \Bigr\}.
\]
Then there exist constants $C(l,\h)$, $0\le l\le d$, such that, for all $K>0$,
\eqs
   \pr[\{E_{1l}^\h(v;K)\}^c \giv \tH(v) ]
        &\le& C(l,\h)K^{-2} W(v) e^{-\l(1-2\h) v}.
\ens
\end{lemma}

\proof
Combining \Ref{ADB-Wj-repn} with~\Ref{ADB-tH-is-Markov}, it follows that 
$\law\bigl((W_0(s),\ldots,W_d(s)),\, s\ge v \giv \hff_v\bigr)$ depends on~$\hff_v$ only through the value of~$\tH(v)$.  Then,
noting that, for $r+\h \le 1$, $1\le l\le d$ and for any $w > t\ge v$,
\[
    \sup_{t \le  s \le w} \{e^{-\l s(1-r-\h)} |W_l(s) - W_l(v)| \}
         \Le  e^{-\l t(1-r-\h)} \sup_{t \le  s \le  w} |W_l(s) - W_l(v)|,
\]
and using Kolmogorov's inequality on the real and imaginary parts of~$W_l$, it follows 
that
\eqs
   \lefteqn{\pr\Bigl[\sup_{t \le  s \le w}\{ e^{-\l s(1-r_l-\h)} |W_l(s) - W_l(v)|\} \ge K 
                    \,\Big|\, \tH(v) \Bigr]} \\
   &&\quad\Le 4K^{-2} e^{-2\l t(1-r_l-\h)} \ex\{|W_l(w) - W_l(v)|^2 \giv \tH(v) \}\,.
\ens
For $r_l > 1/2$, taking $w=\infty$, it follows from~\Ref{AB-MG-variance} that
\eqs
   \lefteqn{\pr\Bigl[\sup_{s \ge v}\{ e^{-\l s(1-r_l-\h)} |W_l(s) - W_l(v)|\} \ge K 
                    \,\Big|\, \tH(v) \Bigr]} \\
   &&\Le 4K^{-2} e^{-2\l v(1-r_l-\h)} W(v) e^{-\l v(2r_l-1)}/(2r_l-1) 
   \Eq 4K^{-2}W(v) e^{-\l v(1-2\h)}/(2r_l-1).
\ens
For $r_l = 1/2$, taking $t = v + j\l^{-1}$ and $w = v + (j+1)\l^{-1}$, it follows from~\Ref{AB-MG-variance} that
\eqs
   \pr\Bigl[\sup_{t \le s \le w}\{ e^{-\l s(1-r_l-\h)} |W_l(s) - W_l(v)|\} \ge K 
                    \,\Big|\, \tH(v) \Bigr]  \Le 4K^{-2} W(v) e^{-(\l v + j)(1-2\h)}(j+1) ,
\ens
and adding over $j \in \Z_+$ gives
\eqs
   \pr\Bigl[\sup_{s \ge v}\{ e^{-\l s(1-r_l-\h)} |W_l(s) - W_l(v)|\} \ge K 
                    \,\Big|\, \tH(v) \Bigr]
   \Le  \frac{4W(v) e^{-\l v(1-2\h)}}{K^2(1 - e^{-(1-2\h)})^2}\,.
\ens
For $r_l < 1/2$, taking $t = v + j\l^{-1}$ and $w = v + (j+1)\l^{-1}$, it follows from~\Ref{AB-MG-variance} that
\eqs
   \pr\Bigl[\sup_{t \le s \le w}\{ e^{-\l s(1-r_l-\h)} |W_l(s) - W_l(v)|\} \ge K 
                    \,\Big|\, \tH(v) \Bigr]
   \Le  \frac{4W(v)e^{-(\l v + j)(1-2\h)}\, e^{1 - 2r_l}}{K^2(1 - 2r_l)}\,,
\ens
and adding over $j \in \Z_+$ gives
\eqs
   \pr\Bigl[\sup_{s \ge v}\{ e^{-\l s(1-r_l-\h)} |W_l(s) - W_l(v)|\} \ge K 
                    \,\Big|\, \tH(v) \Bigr]
   \Le  \frac{4e W(v)e^{-\l v(1-2\h)}}{K^2(1 - e^{-(1-2\h)})(1 - 2r_l)}.
\ens

For $l=0$, the result is proved in analogous fashion, starting from
\[
   \sup_{t \le s \le t + \l^{-1}} \{ e^{\l\h s}|W(s) - W(\infty)|\}
    \Le 2e^{\h(\l t + 1)}\sup_{s\ge t}|W(s) - W(t)|,
\]
and observing that, from~~\Ref{AB-MG-variance},
\eqs
  \phantom{X} \pr\bigl[ \sup_{s\ge t}|W(s) - W(t)| > a \giv \tH(v) \bigr] \Le a^{-2}\ex\{W(t)e^{-\l t} \giv \tH(v)\}
          \Eq a^{-2}e^{-\l t}W(v).   \phantom{X}\halmos
\ens


\medskip
As a result of this lemma, we can sharpen~\Ref{ADB-H-asymp} by giving an explicit bound on the error
made when approximating $e^{-\l t}H_j(t)$ by $W(v)/(d+1)$ for any $t\ge v$. To state the bound,
we define
\eq\label{AB-C-def}
  Q(v) \Def d+2+\sum_{l=1}^d e^{-\l(1- r_l-\h) v}\, |W_l(v)|;\qquad  E_1^\h(v) \Def \bigcap_{l=0}^d E_{1l}^\h(v;1),
\en
noting that, on $E_1^\h(v)$, $Q(t) \le Q(v) + d$ for all $t \ge v$.
Then for all $t\ge v$ and $0\le j\le d$, and if $\h < \rrr(d)$, we have
\eqa
   \lefteqn{\left|e^{-\l t}H_j(t) - \frac{W(v)}{(d+1)}\right|} \non\\
   &\le& \frac1{d+1} \Blb |W(t) - W(v)| + \sum_{l=1}^d e^{-\l(1- r_l) t}\, \{|W_l(v)| + |W_l(t) - W_l(v)|\}\Brb
              \non\\
   &\le& \frac1{d+1} \Blb \sum_{l=1}^d e^{-\l(1- r_l) t}\, |W_l(v)| +
              (d+2) e^{-\l\h t} \Brb \Le \frac{e^{-\l\h v}Q(v)}{d+1}\,, \label{AB-H-as-W-2}
\ena
on~$E_1^\h(v)$. Furthermore, from Lemma~\ref{AB-W-bounds},
\eq\label{AB-H-as-W-prob}
    \pr[\{E_1^\h(v)\}^c \giv \tH(v)] \Le \th_1(v) \Def W(v) e^{-\l(1-2\h) v} 
                  \Bigl( C(0,\h) + \sum_{l=1}^d C(l,\h)\Bigr).
\en

\subsection{Approximating an integral representation of $W(v+t) - W(v)$}
The aim of this section is to prove an approximation theorem, when~$v$ is large, for the process 
$X_v\uo(t) := W(v+t) - W(v)$ in $t\ge0$.
We recall \Ref{AB-incidence-vol} and~\Ref{AB-kappa-vol}, and use the representation~\Ref{AB-Wj-def}, writing
\eqa
    X_v\uo(t) &=& \int_v^{v+t} e^{-\l u}\{M_0(du) - H_{d+1}(du)\} \, \non\\
        &=&   \int_{H_{d+1}(v)}^{H_{d+1}(v+t)} e^{-\l H_{d+1}^{-1}(w)} \{Z\ui(dw) - dw\}\,,
       \label{AB-W-difference}
\ena
where $Z\ui$ is a unit rate Poisson process, with increments independent of~$\hff_v$, starting with 
$Z\ui(H_{d+1}(v)) = M_0(v) = H_0(v)$,
and where~$H_l(u)$, $l\ge0$, are constructed in $u\ge v$ from the Poisson process~$Z\ui$, using 
\Ref{AB-H-eqns} and~\Ref{AB-kappa-vol}, with initial values $H_l(v)$, $0\le l\le d$. Once again, the 
process~$X_v\uo$ depends on its past~$\hff_v$ only through~$\tH(v)$. Since the expression~\Ref{AB-W-difference} 
is too complicated to use directly,  we simplify it in a series of stages.

We start by approximating~$H_{d+1}^{-1}(w)$ in $w \ge H_{d+1}(v)$.  In view of~\Ref{AB-H-as-W-2},
we have $H_{d+1}(t) \approx e^{\l t}W(v)/(d+1)$, or $w \approx e^{\l H_{d+1}^{-1}(w)}W(v)/(d+1)$;
the precise result is as follows.   Note that, for our purposes, $\g^\h(v)$ can be thought of as small.

\begin{lemma}\label{AB-H-inverse}
Fix any $\h < \rrr(d)$.  Then, on the event~$E_1^\h(v)$, we have
\[
   \frac{W(v)(1-\g^\h(v))}{w(d+1)} \Le e^{-\l H_{d+1}^{-1}(w + H^*(v))}
     \Le \frac{W(v)(1+\g^\h(v))}{w(d+1)},
\]
for all $w \ge \{W(v)/(d+1)\}e^{\l v}$, where $\g^\h(v) := (d+1)\{Q(v)/W(v)\}e^{-\l\h v}$,
$H^*(v) := H_{d+1}(v) - e^{\l v}W(v)/(d+1)$, and $Q(v)$ is as defined in~\Ref{AB-C-def}.
\end{lemma}

\proof
We begin by noting that $H_{d+1}(u) = \int_0^u  \l H_d(t)\,dt$, so that, from~\Ref{AB-H-as-W-2}, 
for $u\ge v$,
\eqa\label{AB-H-bounded}
 \lefteqn{ \left|H_{d+1}(u) - H_{d+1}(v) - (e^{\l (u-v)} - 1)e^{\l v}\frac{W(v)}{d+1}\right| }\non\\
     &&\qquad\Le \int_v^u \l e^{\l t}
       \Blb  \sum_{l=1}^d |W_l(v)|e^{-\l(1-r_l)t} + (d+2)e^{-\l\h v} \Brb \,dt\non\\
    &&\qquad\Le Q(v) e^{\l (u-v)} e^{\l(1-\h)v}\,.            
\ena
So, defining 
\eq\label{AB-t(s)-def}
       t_v(s) \Def \l^{-1}\log\Blb 1 + \frac{s(d+1)}{e^{\l v}W(v)} \Brb \quad\mbox{and}\quad
     t_v^{-1}(u) \Def \frac{e^{\l v}W(v)}{d+1}\,(e^{\l u} - 1) ,
\en     
it follows that, on $E^\h_1(v)$, 
\eqa
    \lefteqn{|\{H_{d+1}(t_v(s)+v) - H_{d+1}(v)\} - s|}\non\\ [1ex]
       &&\Le Q(v)e^{\l(1-\h)v}\Blb 1 + \frac{s(d+1)}{e^{\l v}W(v)} \Brb \ =:\ h_v(s). 
                    \label{AB-delta-H-bnd-1}
\ena
Now substitute $s = t_v^{-1}(u)$ into~\Ref{AB-delta-H-bnd-1}
for $u\ge0$, giving
\[
   \frac{W(v)}{d+1}\,e^{\l (u+v)}(1-\g^\h(v)) + H^*(v) \Le H_{d+1}(u+v) 
     \Le \frac{W(v)}{d+1}\,e^{\l (u+v)}(1+\g^\h(v)) + H^*(v).
\]
Writing $w = H_{d+1}(u+v)$ and inverting, it then follows immediately that
\[
    \l^{-1}\log\Blb \frac{(w-H^*(v))(d+1)}{W(v)(1+\g^\h(v))} \Brb \Le H_{d+1}^{-1}(w)
      \Le \l^{-1}\log\Blb \frac{(w-H^*(v))(d+1)}{W(v)(1-\g^\h(v))} \Brb,
\]
establishing the lemma.
\ep
\ignore{
\[
   \frac{W(v)(1-\g(v))}{w(d+1)} \Le e^{-\l H_{d+1}^{-1}(w + H^*(v))}
     \Le \frac{W(v)(1+\g(v))}{w(d+1)}.
\]
}

This now allows~\Ref{AB-W-difference} to be rewritten in the form
\eq\label{AB-W-difference-2}
   X_v\uo(t) \Eq \int_0^{H_{d+1}(v+t) - H_{d+1}(v)} 
         e^{-\l H_{d+1}^{-1}(w + H_{d+1}(v))}(Z\ut(dw) - dw),
\en
where~$Z\ut$ is a unit rate Poisson process, with respect to which both upper limit and integrand 
are predictable, the latter being decreasing in~$w$ and bounded between 
\eq\label{AB-W-integrand-bnd}
    \frac{W(v)(1-\g^\h(v))}{w(d+1) + W(v)e^{\l v}} \quad\mbox{and}\quad 
      \frac{W(v)(1+\g^\h(v))}{w(d+1) + W(v)e^{\l v}}\,,
\en
for all $w\ge0$,
on the event~$E_1^\h(v)$. In order to show that we can replace both the integrand and the upper limit of 
integration in~\Ref{AB-W-difference-2} with simpler expressions, without making too great an error, we use 
Lemma~\ref{AB-process-bounds} from the Appendix.

We first replace the integrand in~\Ref{AB-W-difference-2}, showing 
that~$X_v\uo$ is close to~$X_v\ui$, defined by
\eq\label{AB-X1-def}
   X_v\ui(t) \Def \int_0^{H_{d+1}(v+t) - H_{d+1}(v)} \frac{W(v)}{w(d+1) + W(v)e^{\l v}}\,(Z\ut(dw) - dw),
\en
using~\Ref{AB-W-integrand-bnd}.  We set
\[
    v_-(\h) \Def \max\bigl\{0,[\l(1-\h)]^{-1}\log\{e^{-2}(d+1)\}\bigr\}.
\]

\begin{lemma}\label{AB-smooth-integrand}
With the above definitions, for any $\h < \rrr(d)$ and any $v \ge v_-(\h)$,
we have
\eqs
 \lefteqn{
     \pr\Bigl[ e^{\l v/2} \sup_{t \ge 0}|X_v\uo(t) - X_v\ui(t)| > \{W(v)Q(v)\g^\h(v)\}^{1/2} 
        \giv \tH(v) \Bigr]    } \\ 
   && \Le \th_2(v) \Def \th_1(v) + \tilde\th_2(v), \phantom{XXXXXXXXXXXXXXXX}
\ens
where $\th_1(v)$ is as in~\Ref{AB-H-as-W-prob}, and $\tilde\th_2(v) := 2e^{-W(v)e^{\l\h v}/\{2e\}}$.
\end{lemma} 
 
\proof
It follows from~\Ref{AB-W-difference-2} that $X(t) = X_v\uo(t) - X_v\ui(t)$ is an integral of the form
considered in Lemma~\ref{AB-process-bounds}, albeit with a random upper limit, and its corresponding function~$F$
satisfies
\eq\label{AB-F-bound}
     |F(u)| \Le G(u) \Def \frac{\g^\h(v)W(v)}{u(d+1) + W(v)e^{\l v}},\quad\mbox{for all}\ u\ge0,
\en
on $E_1^\h(v)$, in view of~\Ref{AB-W-integrand-bnd}.  We can thus apply Lemma~\ref{AB-process-bounds}
to the process~$\tX$ with $\tF(t) := F(t)\bone\{|F(u)| \le G(u),\,0\le u < t\}$ and with~$\tG(u) := G(u)$ as
in~\Ref{AB-F-bound}, noting that then, recalling~\Ref{AB-H-as-W-prob},
\[
    \pr[X(t) = \tX(t)\ \mbox{for all}\ t\ge0 \giv \hff_v] \ \ge\ \pr[E_1^\h(v) \giv \tH(v)] \ \ge\ 1 - \th_1(v). 
\]
Now, from~\Ref{AB-F-bound}, we have $\tG_2(0,\infty) = \{\g^\h(v)\}^2 \{W(v)/(d+1)\} e^{-\l v}$.  We can
then choose $a := e^{-\l v/2}\{W(v)Q(v)\g^\h(v)\}^{1/2}$ in Lemma~\ref{AB-process-bounds},  because 
\[
   a \Le e\tG_2(0,\infty)/\tG^*(0,\infty) \Eq e\g^\h(v)\{W(v)/(d+1)\}
\]
if $v \ge v_-(\h)$, and the result follows.
 \ep

The next step is to simplify the upper limit in~\Ref{AB-X1-def}, using Lemma~\ref{AB-process-bounds} to
show that, with $t_v(s)$ as defined in~\Ref{AB-t(s)-def},
 $(X_v\ui(t_v(s)),\,s\ge0)$ is close to the process~$(X_v\ut(s),\,s\ge0)$ given by
\eq\label{AB-X2-def}
    X_v\ut(s) \Def \int_0^s \frac{W(v)}{w(d+1) + W(v)e^{\l v}}\,(Z\ut(dw) - dw).
\en
 For this, we need to control
$\sup_{s\ge0,\,|z|<h_v(s)} |X_v\ut(s+z) - X_v\ut(s)|$, for~$h_v(s)$ defined 
in~\Ref{AB-delta-H-bnd-1}.

\begin{lemma}\label{AB-X2-smooth}
With the definitions given in \Ref{AB-delta-H-bnd-1}, \Ref{AB-X1-def} and \Ref{AB-X2-def}, 
and for any $\h < \rrr(d)$, we have
\eqs
   \lefteqn{ \pr\Bigl[ e^{\l v/2}\sup_{s\ge0,\,|z|<h_v(s)} |X_v\ut(s+z) - X_v\ut(s)| > 4\e^\h(v) 
               \,\Big|\, \tH(v) \Bigr] I[E_{21}^\h(v)]}\\
    && \qquad\Le \th_3(v) \Def  \Blb 2\left[1 + \frac{W(v)}{g(v)}\right] 
     + \frac{8e\,e^{2\l\h v/3}}{Q(v)(d+1)^{2}}\Brb   e^{-W(v)e^{\l\h v/3}/\{2e(d+2)\}},
\ens
where $\e^\h(v) := \{W(v)Q(v)\}^{1/2}e^{-\l\h v/3}$, $g(v) := Q(v)(d+2)e^{-\l\h v}$ and 
\eq\label{AB-E21-def}
   E_{21}^\h(v) \Def \{W(v) \le e^2(d+2)^2 Q(v)e^{\l v/3}\}\cap\{Q(v) \le 2e(d+1)^{-2}e^{2\l\h v/3}\}
    \ \in\ \s(\tH(v)).
\en
\end{lemma}

\proof
We consider the ranges $0 \le s \le W(v)e^{\l v}$ and $s >  W(v)e^{\l v}$ separately.
In the first range of~$s$, define $s_j := j e^{\l v} g(v)$
for $0 \le j \le M := \lfloor W(v)/g(v) \rfloor$, and set $s_{M+1} := W(v)e^{\l v}$:
then $s_{j+1} - s_j \ge h_v(s_j)$ for each~$j$.
By Lemma~\ref{AB-process-bounds}, with $G(u)$ the constant~$e^{-\l v}$ and $a := e^{-\l v/2}\e^\h(v)$, we have
\[
   \pr\Bigl[\sup_{s_j \le s \le s_{j+1}} e^{\l v/2}|X_v\ut(s) - X_v\ut(s_j)| > \e^\h(v) 
         \giv \tH(v) \Bigr]I[E_{21}^\h(v)] \Le 2\exp\{-\e^\h(v)^2 / (2eg(v))\}, 
\]
for $0 \le j \le M$, since $a  \le eg(v) = eG(s_j)(s_{j+1}-s_j)$ 
on~$E_{21}^\h(v)$.  Hence,
by a standard argument,
\eqa
 && \pr\Bigl[\sup_{0\le s \le W(v)e^{\l v},\,|z| < h_v(s)} e^{\l v/2}|X_v\ut(s+z) - X_v\ut(s)| > 3\e^\h(v)
                \,\Big|\, \tH(v) \Bigr]I[E_{21}^\h(v)]\non\\
    &&\qquad\qquad \Le 2\{1 + W(v)/g(v)\} \exp\{-W(v)e^{\l\h v/3}/\{2e(d+2)\}\}.\label{AB-range-1}
\ena
In the second range of~$s$, we define
\[
   s_j \Def W(v) e^{\l v}(1 + \tg(v))^j, \quad \mbox{where}\quad \tg(v) \Def g(v)(d+1)/W(v),
\]
noting that $s_{j+1}-s_j = s_j\tg(v) \ge h_v(s_j)$.
By Lemma~\ref{AB-process-bounds} with $G(u) := s_j^{-1}\{W(v)/(d+1)\}$, we have
\eqs
  && \pr\Bigl[\sup_{s_j \le s \le s_{j+1}} e^{\l v/2}|X_v\ut(s) - X_v\ut(s_j)| > \e^\h(v) 
                \,\Big|\, \tH(v) \Bigr] I[E_{21}^\h(v)]\\
   &&\qquad\qquad \Le 2\exp\{-\{\e^\h(v)\}^2 (d+1)^2(1 + \tg(v))^j/ (2eW(v)\tg(v))\}, \quad j\ge0,
\ens
since $a := e^{-\l v/2}\e^\h(v) \le eg(v) = e\{W(v)/(d+1)\}\tg(v) = eG(s_j)(s_{j+1}-s_j)$ on~$E_{21}^\h(v)$,
and hence 
\eqa
  &&\pr\Bigl[\sup_{s \ge W(v) e^{\l v},\,|z| < h_v(s)} e^{\l v/2}|X_v\ut(s+z) - X_v\ut(s)| > 4\e^\h(v) 
                \,\Big|\, \tH(v) \Bigr] I[E_{21}^\h(v)]\non\\
    &&\ \Le 2 \exp\{-W(v)(d+1) e^{\l\h v/3}/ \{2e(d+2)\}\} 
                    \sum_{j\ge0} \exp\{-j\{\e^\h(v)\}^2(d+1)^2/(2eW(v))\} \non\\
    &&\ \Le \frac{8e\,e^{2\l\h v/3}}{Q(v)(d+1)^{2}} \exp\{-W(v)(d+1) e^{\l\h v/3}/ \{2e(d+2)\}\}, 
           \label{AB-range-2}
\ena
since also $\{\e^\h(v)^2\}(d+1)^2/(2eW(v)) \le 1$ on~$E_{21}^\h(v)$.  We need $4\e^\h(v)$ here as the bound on the 
supremum difference, rather than the usual $3\e^\h(v)$, because
it is possible to have $s(1 - \tg(v)) < s_{j-1}$ for some $s_j< s < s_{j+1}$;  however, it then has to be 
the case that, for such~$s$, $s(1-\tg(v)) \ge s_{j-2}$ if $\tg(v) \le 1/2$, which is
the case on~$E_{21}^\h(v)$.
\ep

In view of Lemma~\ref{AB-X2-smooth} and~\Ref{AB-delta-H-bnd-1}, we immediately have the following corollary.

\begin{corollary}\label{AB-X1X2-cor}
With the definitions of Lemma~\ref{AB-X2-smooth},
$$
 \pr\Bigl[ e^{\l v/2}\sup_{s\ge0} |X_v\ui(t_v(s)) - X_v\ut(s)| > 4\e^\h(v) 
               \,\Big|\, \tH(v) \Bigr] I[E_{21}^\h(v)] \Le \th_1(v) + \th_3(v).
$$
\end{corollary}

We now show that $X_v\ut$ is close in distribution to the process~$X_v\uh$ defined by
\eq\label{AB-X3-def}
   X_v\uh(s) \Def \int_0^s \frac{W(v)}{w(d+1) + W(v)e^{\l v}}\,B(dw),
\en
where, for the integrator, the compensated Poisson process $Z\ut(w) - w$ from~$X_v\ut$ has been replaced by a 
standard Brownian motion~$B(w)$.  Note that~$e^{\l v/2}X_v\uh$ is itself just a time-changed Brownian motion: 
\eq\label{AB-brownian-repn}
   \Bl\{(d+1)/W(v)\}^{1/2}\,X_v\uh\bigl(\{W(v)/(d+1)\}e^{\l v}s\bigr),\,s\ge0\Br 
             \ =_d\ \bigl(B(s/(s+1)),\,s\ge0\bigr), 
\en
and so, conditional on~$W(v)$, $X_v\uh(\infty) \sim \nn(0,W(v)/(d+1))$.

\begin{lemma}\label{AB-KMT}
Fix $r\ge1$.  Then there are constants $c_{r1}$ and~$c_{r2}$, depending only on~$d$, with the following properties.
For all~$v$ such that $\l v \ge c_{r1}$,
it is possible to construct $X_v\ut$ and~$X_v\uh$ on the same probability space, in such
a way that
\[
   \pr\Bigl[e^{\l v/2}\sup_{s\ge0}|X_v\uh(s) - X_v\ut(s)| \ge c_{r2}(1+ W(v)) \l v e^{-\l v/2} \Bigr] \Le 
           \th_4(v) \Def e^{-3r\l v}.
\]
\end{lemma}

\proof
For any $r \ge 1$, there are constants $C_r,K_r$ with the property that, for any $n\ge1$, 
a standard Poisson process~$Z$ and
a standard Brownian motion~$B$ can be constructed on the same probability space in such a way that
$\pr[A_r^c(n)] \le  K_r n^{-(r+1)}$, where
\[
   A_r(n) \Def \left\{ \sup_{0\le s\le n} \frac{|Z(s) - s - B(s)|}{\log n} \le C_r \right\}.
\]
This follows from Koml\'os, Major \& Tusn\'ady~(1975, Theorem~1 (ii)), together with elementary exponential bounds 
for the fluctuations of the standard Poisson process and Brownian motion over the time interval $[0,1]$. Fix~$r$, 
and take $n := e^{3\l v}$ for $v \ge v_1$, where~$v_1$ is chosen so that $e^{3\l v_1} \ge 2K_r$, implying that 
$\pr[A_r^c(n)] \le  \frac12 e^{-3r\l v}$. Then use the corresponding choices of $Z$ and~$B$ to realize 
$X_v\ut$ and~$X_v\uh$, which we express, by partial integration, in the form 
\eqa
  X_v\ut(s) &=& W(v) \Blb \frac{Z(s)-s}{s(d+1) + W(v)e^{\l v}} 
                         + \int_0^s \frac{Z(u)-u}{(u(d+1) + W(v)e^{\l v})^2}\,du \Brb, \non\\
  X_v\uh(s) &=& W(v) \Blb \frac{B(s)}{s(d+1) + W(v)e^{\l v}} 
                         + \int_0^s \frac{B(u)}{(u(d+1) + W(v)e^{\l v})^2}\,du \Brb.
   \label{AB-X-partial-form}
\ena
Taking the difference, it is immediate that, for $0 \le s \le e^{3\l v}$ and on~$A_r(e^{3\l v})$,
\eqs
   \frac{|Z(s)-s-B(s)|}{s(d+1) + W(v)e^{\l v}} &\le& C_r\, \frac{3\l v}{ W(v)e^{\l v}} 
\ens
and that
\[
    \int_0^s \frac{|Z(u)-u-B(u)|}{(u(d+1) + W(v)e^{\l v})^2}\,du \Le 
      C_r \, \frac{3 \l v}{W(v)(d+1) e^{\l v}} \,.
\]
This shows that, on~$A_r(e^{3\l v})$,
\[
    e^{\l v/2}|X_v\uh(s) - X_v\ut(s)| \Le  6C_r \l v e^{-\l v/2} 
            \quad\mbox{for}\quad   0 \le s \le e^{3\l v}.
\]
Then, taking $F(u) = W(v)/\{u(d+1)+W(v)e^{\l v}\}$, $a = eC_r \{W(v)/(d+1)\} \l v e^{-\l v}$,  $t_1=e^{3\l v}$
and $t_2=\infty$
 in Lemma~\ref{AB-process-bounds}, with the choice of~$a$ permissible
for all $v\ge v_2$, where~$v_2 \ge \l^{-1}$ is chosen such that $\l v_2 e^{-\l v_2} \le 1/C_r$,
we have
\[
  \pr\Bigl[ \sup_{e^{3\l v} \le s < \infty}|X_v\ut(s) - X_v\ut(e^{3\l v})| 
                                            > eC_r \{W(v)/(d+1)\} \l v e^{-\l v} \Bigr]
      \Le 2\exp\{-(e/2)(C_r \l v)^2 e^{\l v}\}.
\]
The same bound is satisfied also for $\sup_{e^{3\l v} \le s < \infty}|X_v\uh(s) - X_v\uh(e^{3\l v})|$,
as can be deduced from the representation~\Ref{AB-brownian-repn}.
Now choose~$v_3 \ge \l^{-1}$ so that $8\exp\{-(e/2)(C_r\l v_3)^2 e^{\l v_3}\} \le e^{-3r\l v}$, 
and set $v_0 :=\max\{v_1,v_2,v_3\}$.
\ep

Summarizing the conclusions Lemmas \ref{AB-smooth-integrand} and~\ref{AB-KMT} and of 
Corollary~\ref{AB-X1X2-cor}, we have the following theorem.  In the error terms,
$\th_1(v)$ is defined in~\Ref{AB-H-as-W-prob}, $\th_2(v)$ in Lemma~\ref{AB-smooth-integrand},
$\th_3(v)$ in Lemma~\ref{AB-X2-smooth} and $\th_4(v)$ in Lemma~\ref{AB-KMT}.

\begin{theorem}\label{AB-branching-appxn}
With the definitions \Ref{AB-Wj-def}, \Ref{AB-t(s)-def} and~\Ref{AB-X3-def}, fixing
any $\h < \rrr(d)$, we can construct $W$ and a time changed Brownian motion~$X_v\uh$
on the same probability space, in such a way that, for all $v \ge \l^{-1}c_{1*}$,
\[ 
   \pr\Bigl[e^{\l v/2}\sup_{u\ge0}|\{W(u+v) - W(v)\} - X_v\uh(t_v^{-1}(u))| > K(v)e^{-\l\h v/3}
             \,\Big|\, \tH(v)  \Bigr] I[E_{21}^\h(v)] \Le \sum_{i=1}^4 \th_i(v),
\]
where $K(v) := 4\{W(v)Q(v)\}^{1/2} + Q(v)\sqrt{d+1} + c_{2*}(1+W(v)e^{-\l v/3})$, 
$E_{21}^\h(v) \in \s(\tH(v))$ is as defined in~\Ref{AB-E21-def}, and the
constants $c_{1*}$ and~$c_{2*}$, which depend only on~$d$, can be deduced from Lemma~\ref{AB-KMT} with~$r=1$.
\end{theorem}

\subsection{Consequences for the gossip process}
Theorem~\ref{AB-branching-appxn} is not yet in a form easily applied to the gossip process.
To start with,  
the statement of the theorem involves the $\s(\tH(v))$-measurable random variables $W(v)$, $Q(v)$, $K(v)$
and $\th_i(v)$, $1\le i\le 4$, and it is useful to have some idea of their magnitude.  It is also useful
to specify how big the probability~$\pr[E_{21}^\h(v)]$ may be.
To derive appropriate statements, we begin with the random elements $W(v)$ and~$W_l(v)$, $1\le l\le d$.

\begin{lemma}\label{AB-W-Wj-control}
For any $0 < \h < \rrr(d)$, we have
\eq\label{AB-Wj-bounds-1}
   \pr[e^{-\l(1- r_l-\h) v}\, |W_l(v)| > 2] \Le \left\{
    \begin{array}{ll}
          e^{-2\l v(1- r_l-\h) } (2r_l-1)^{-1} &\quad\mbox{if}\quad r_l > 1/2; \\
          \l v\, e^{-\l v(1-2\h)} &\quad\mbox{if}\quad r_l = 1/2; \\
          e^{-\l v(1-2\h)}(1-2r_l)^{-1} &\quad\mbox{if}\quad r_l < 1/2,
    \end{array} \right.
\en
for $1\le l\le d$.  Furthermore, for any $s > 0$,
\eq\label{AB-W0-bounds}
    \pr[W(v) \ge 1 + s] \Le  s^{-2} \quad \mbox{and}\quad
    \pr[W(v) \le s] \Le \exp\Blb - \frac{\{\log_+(w_0/s)\}^{d+1}}{2\,(d+1)!} \Brb,
\en
for a suitably chosen $w_0 > 0$.
\end{lemma}
 
\proof
The first part follows from~\Ref{AB-MG-variance} and Chebyshev's inequality,
and, for~$W(v)$, the bound on the upper tail holds because 
\hbox{$\var W(v) \le \var W(\infty) \le 1$} and $\ex W(v) = 1$.  For the lower tail,
note that $W(\infty) > 0$ a.s., so that, because~$W(\cdot)$ is c\`adl\`ag and
positive on~$\re_+$, we have $W_* := \inf_{t > 0}W(t) > 0$ a.s.\ also.  Suppose that
$w_0 > 0$ is such that $\pr[W_* \ge w_0] \ge 1/2$.  Then, for $0 < x \le w_0$, $W(t) > x$ if any of the offspring
of the initial individual that are born before time~$t_x$ generate families with~$W_* > w_0$,
where $e^{-\l t_x} = x/w_0$.  The probability that there are no such offspring is just
$\exp\{- \r \vva  t_x^{d+1}/\{2(d+1)\}\}$.  Hence, for $t\ge t_x$ and $x \le w_0$,
\[
   \pr[W(t) \le x] \Le \exp\Blb -\frac{\r \vva  \{\log(w_0/x)\}^{d+1}}{2\l^{(d+1)}(d+1)} \Brb
                   \Eq \exp\Blb - \frac{\{\log(w_0/x)\}^{d+1}}{2\,(d+1)!} \Brb.
\]
\ep

In view of~\Ref{AB-C-def}, if $0 < \h < \rrr(d)$, then $Q(v) \le 3(d+1)$
on the event
\eq\label{AB-E22-def}
  E_{22}^\h(v) \Def \bigcap_{l=1}^d \{e^{-\l(1- r_l-\h) v}\, |W_l(v)| \le 2\},
\en 
and the first part~\Ref{AB-Wj-bounds-1} of Lemma~\ref{AB-W-Wj-control} directly implies that
\eq\label{AB-E22-bnd}
    \pr[\{E_{22}^\h(v)\}^c] \Le  c(d)(1 + \l v\bone_{\{d=6\}})\,e^{-2\l v(\rrr(d)-\h)},
\en
for a suitable constant~$c(d)$; of course, by definition, $Q(v) \ge d+2$.  
The second part of Lemma~\ref{AB-W-Wj-control} implies that $E_{23}(v) := \{W(v) \le 1 + e^{\l\h v/3}\}$
is such that $\pr[\{E_{23}(v)\}^c] \le e^{-2\l\h v/3}$.  From these observations and~\Ref{AB-E21-def}, 
it follows that
\[
     E_{22}^\h(v) \cap E_{23}(v) \ \subset\  E_{21}^\h(v),
\]
if $v$ is such that $e^{2\l\h v/3} \ge (d+1)^3$, and hence, for such~$v$,
\eq\label{AB-E21-bnd}
    \pr[\{E_{21}^\h(v)\}^c] \Le c(d)(1 + \l v) e^{-2\l v(\rrr(d)-\h)} + e^{-2\l\h v/3};
\en
in addition, 
$$
          K(v)e^{-\l\h v/3} \Le \sqrt{d+1}\,\{4\sqrt2  + 2(d+1) + 3c_{2*}\}\,e^{-\l\h v/6}
$$
on $E_{22}^\h(v) \cap E_{23}(v)$ also.  

For the quantities~$\th_i$, $1\le i\le 4$, note that, from~\Ref{AB-H-as-W-prob}, 
\eq\label{theta-1-bnd}
   \th_1(v) \Le C(d,\h)e^{-\l v(\rrr(d)-\h)} \quad \mbox{on the event}\ 
            E_{24}^\h(v) \Def \{W(v) \le 1 + e^{\l v(\rrr(d)-\h)}\}, 
\en
and that $\pr[\{E_{24}^\h(v)\}^c] \le e^{-2\l v(\rrr(d)-\h)}$.  Then, as in Lemma~\ref{AB-KMT},
$\th_4(v) = e^{-3\l v}$ if we take $r=1$.  From Lemma~\ref{AB-smooth-integrand}, $\th_2(v) = \th_1(v) + \tth_2(v)$, and
both $\tth_2(v)$ and~$\th_3(v)$, defined in Lemma~\ref{AB-X2-smooth}, are super-exponentially small in~$\l\h v$ on
the event $E_{25}^\h(v) := \{W(v) \ge e^{-\l\h v/6}\}$. 
Finally, by the last inequality in Lemma~\ref{AB-W-Wj-control},
\[
   \pr[\{E_{25}^\h(v)\}^c] \Le \exp\Blb - (1/2)\{\log(w_0) + \l\h v/6\}^{d+1}/(d+1)! \Brb,
\]
which is also super-exponentially small in~$\l\h v$.  Hence, taking
\eq\label{AB-E-final-def}
  E^\h(v) \Def  E_{22}^\h(v) \cap E_{23}(v) \cap E_{24}^\h(v) \cap E_{25}^\h(v),
\en
for which $\pr[\{E^\h(v)\}^c] \le C(d)(\l v e^{-2\l v(\rrr(d)-\h)} + e^{-2\l\h v/3})$,
and assuming that~$v$ is such that $e^{2\l\h v/3} \ge (d+1)^3$, we have the following consequence of
Theorem~\ref{AB-branching-appxn}.  To state it, and for future use,  we define
\eq\label{AB-tlm-def}
    \tlm \Def (3/2)\l^{-1}\log\L,
\en
an upper bound for the  times to be considered in proving the central limit theorem.

\begin{corollary}\label{AB-CLT-cor}
For any $0 < \h < \rrr(d)$ and $v \le \tlm$ such that $e^{2\l\h v/3} \ge (d+1)^3$ and $\l v \ge c_{1*}$, 
there are constants~$C = C(d,\h)$ and~$C' = C'(d)$  and an event $E^\h(v) \in \s(\tH(v))$, with
$\pr[\{E^\h(v)\}^c] \le C' \l v e^{-2\l v(\rrr(d)-\h)} + e^{-2\l\h v/3})$,
such that, for any $u\ge0$ such that $t_{\L}(u) \le \tlm$,
\eqa
   \lefteqn{|\ex\bigl\{f(e^{\l v/2}\{W(u+v) - W(v)\}) \giv \hff_v  \bigr\}
       - \ex\bigl\{ f(e^{\l v/2} X_v\uh(t_v^{-1}(u))) \giv \hff_v  \bigr\}| I[E^\h(v)] } \non\\
   &&\Le C\{e^{-\l\h v/6} + e^{-\l v(\rrr(d)-\h)}\}, \phantom{XXXXXXXXXXXXXXXXXXXXX}
             \label{AB-Lipschitz-approx}
\ena
uniformly for all $f\in\Fbw$.
\end{corollary}

Taking any $c_0,\ldots,c_d \in \re_+$  and setting 
$C(x) := \sum_{l=0}^d c_l x^l$, we also observe from Lemma~\ref{ADB-martingales} that
\eq
   \Bigl| \sum_{l=0}^d c_l e^{-\l s}H_l(s) - (d+1)^{-1}C(1) W(s) \Bigr|
   \Le \frac{C(1)}{d+1} \,\sum_{l=1}^d e^{-\l(1- r_l)s} |W_l(s)| 
    \Le C(1) e^{-\l\h s}        \label{AB-H-sum-approx}
\en
on $E_{22}^\h(s)$, the probability of whose complement is bounded in~\Ref{AB-E22-bnd}.

\section{The central limit theorem}\label{CLT}
In this section, the central limit theorem is proved much as outlined in the introduction. 
With $\s^2_L(v,u) := \var\{L_{t_\L(u)}/L\giv \ff_v\}$, we show in Lemma~\ref{AB-variance-lemma} that
\eq\label{approx-1}
   \ex\Blb \bigl|(L_{t_\L(u)}/L) - \ex\{L_{t_\L(u)}/L \giv \ff_s\}\bigr|\Giv \ff_v\Brb \ \ll\ \s_L(v,u),
\en
if~$s$ is chosen to be sufficiently long after~$v$.  The approximation of $\ex\{L_{t_\L(u)}/L \giv \ff_s\}$ as a Poisson
probability is then accomplished in Lemma~\ref{AB-first-representation}, with an error that is 
small if ${t_\L(u)}-s$ is sufficiently 
large.  Lemmas \ref{AB-M-approx}--\ref{AB-W-to-Wsv} approximate the mean of the Poisson distribution by successively 
simpler quantities, and
bound the errors involved in the approximations.  The combined result of these steps is summarized
in Corollary~\ref{summary-1}, showing that, given~$\ff_v$, the distribution of $L_{t_\L(u)}/L$ is close to
that of $\ell(\log[\hc_d W(s,v)]+ u)$.   

Now the normalized difference $e^{\l v/2}(W(s,v)-W(v,v))$ can be shown, 
using Corollary~\ref{AB-CLT-cor}, to have a normal approximation.  Because of the normalization,
it is important at this point to check that the approximation errors in the previous steps are all much 
smaller than~$e^{-\l v/2}$; this places some restrictions on how large~$v$ may be.
The linearization of the difference  $\ell(\log[\hc_d W(s,v)]+ u) - \ell(\log[\hc_d W(v,v)]+ u)$, needed to show that
it is itself approximately normally distributed, is
accomplished in Lemma~\ref{K-version}, and the final result is given in Theorem~\ref{AB-main-theorem}.

\subsection{Comparisons of processes}
The detailed calculations make heavy use of comparisons between a number of processes, that we justify
in Lemma~\ref{AB-dominating-processes} by realizing them on the same probability spaces. 
The process~$\law$ itself can be realized by starting with the times~$(\bt_j,\,j\ge0)$ of the
branching process~$\bX$, paired with a sequence of independent uniform points $(\bP_j,\,j\ge0)$
of~$\cC$.  This yields a process 
\eq\label{ADB-Y-proc-def}
    Y(t) \Def \{(\bt_j,\bP_j),\,j\in\bJ_t\},\quad t\ge0,
\en 
in terms of which we define 
\eq\label{ADB-bJ-etc-defs}
 \bJ_t \Def \{j \ge 0\colon \bt_j \le t\};\quad \bN_t \Def |\bJ_t|; \quad           
        \bM_t \Def \sum_{j\in\bJ_t} (t-\bt_j)^d.
\en
We can then define the set valued process
\eq\label{ADB-baL-def}
    \baL(t) \Def \bigcup_{j \in \bJ_t} \KK(\bP_j,t-\bt_j),
\en
obtained by taking the unions of the neighbourhoods generated by~$Y(t)$. The process~$Y$
can be augmented to a process~$\tY$ of quadruples, by including a set of pairs~$(K(j),\bQ_j)$, $j\ge0$,
where $0 \le K(j) < j$ and $\bQ_j \in \cC$, 
denoting the subsets from which the long range contacts were made and the positions of the individuals
within them:  given~$Y(\bt_j-)$,
\[
    \pr[K(j)=l] \Eq \frac{\vva_{\bt_j-\bt_l}}{\sum_{l'=0}^j \vva_{\bt_j-\bt_{l'}}},\quad 0\le l < j,
\]
and~$\bQ_j$ is then chosen uniformly from the set $\KK(\bP_{K(j)},\bt_j - \bt_{K(j)})$.
The process~$\law$ is derived from~$\tY$ sequentially, by thinning.   
The pair $(\bt_j,\bP_j)$ is not included in~$\law$ unless
$K(j) = \min\{l \ge 0\colon \bQ_j \in \KK(\bP_{l},\bt_j-\bt_{l})\}$.
This thinning process ensures that, when neighbourhoods overlap in~$\cC$, only
contacts from the neighbourhood that was informed earliest are allowed, ensuring
that the rate of long range transmissions from~$\law_t$ remains equal to~$\r L_t$.
Note that, if $\bP_j \in \law_{\bt_j-}$, the pair $(\bt_j,\bP_j)$ is included in defining~$\law$;
however, it is redundant in~\Ref{ADB-informed-set}, the newly informed individual 
having previously been informed, and it never contributes to further transmission,
because of the definition of the thinning step.  
The resulting set of times and positions we denote by $((\t_j,P_j),\,j\ge0)$, 
with 
\eq\label{ADB-Js-etc-defs}
  J_s \Def \{j\ge0\colon \t_j \le s\};\quad N_s \Def |J_s|;\quad M_s \Def \sum_{j\in J_s}(s-\t_j)^d,
\en
and $\law$ is as given by~\Ref{ADB-informed-set}; it satisfies
$\law_t \subset \baL_t$, with strict inclusion for all large enough times.

The process~$\baL$ acts as a tractable upper bound for~$\law$, and it is
useful also to have tractable lower bounds.  In particular, when calculating the probability that
a neighbourhood~$\KK(P,s)$ intersects~$\law_t$, where~$s$ is fixed and~$P$ is a uniform
random point of~$\cC$, the way in which the neighbourhoods of~$\law_t$ intersect one another
enters in a complicated way.  However, if~$\law_t$ happened to consist of a union
of {\it non-intersecting\/}
neighbourhoods, which were also separated from one another by distance at least~$2s$, then
the probability could be deduced by simply adding the intersection probabilities for the individual
neighbourhoods. Then, because the neighbourhoods~$\KK$ are balls in a geodesic metric space,  
the probability of two neighbourhoods $\KK(P,s)$ and~$\KK(Q,t)$
intersecting, if one or both of $P$ and~$Q$ are chosen uniformly and independently in~$\cC$, is
given by
\eq\label{AB-intersection-prob}
   q_L(s,t) \Eq L^{-1}\vva_{s+t}, 
\en
where~$\vva_{s+t}$ can be estimated in terms of $\vva (s+t)^d$, in view of~\Ref{AB-volume-approx}.
Of course, as~$t$ grows, intersections occur in~$\law_t$, but, at least for a while, their effect may not be
too large.  So the next step is to construct subsets of~$\law_t$ with the necessary separation
properties, and which are amenable to analysis.

Fix any $s,t > 0$, and thin the process~$\tY$ to obtain a set valued process~$\law\ust$
as follows.  Start with $\t_0\ust = 0$ and $P\ust_0 = P_0$, defining 
\[
    \law\ust_u \Def \KK(P_0, u) \quad\mbox{for}\ 0 \le u < \bt_1;
\]
let $R\ust_0 := \emptyset$ denote the initial set of indices of censored points of~$\tY$.  Then
proceed sequentially.  Suppose that the  
quadruples $((\bt_{l},\bP_{l},K(l),\bQ_{l}),\,0\le l\le j-1) \subset \tY$ have already been considered.
If $K(j) \in R\ust_{j-1}$, set $R\ust_j := R\ust_{j-1} \cup \{j\}$ and proceed to the next quadruple;
descendants of censored points are also censored.  If not, 
thin much as in the construction of~$\law$, except that a point~$\bP_j$ is also thinned if it
belongs to $N_{2s+t-\bt_j}(\law\ust_{\bt_j-})$, where, for $V \subset \cC$ 
and $u > 0$,
\eq\label{ADB-LB-thinning-def}
     N_{u}(V) \Def \bigcup_{y \in V} \KK(y,u);
\en
set 
\eq\label{ADB-lawst-rule}
     \law\ust_u \Def \bigcup_{l=0}^j \bone_{\{l \notin R\ust_j\}} \KK(\bP_l,u-\bt_l), \quad \bt_j \le u < \bt_{j+1}.
\en
The extra thinning in~\Ref{ADB-LB-thinning-def}
ensures that the neighbourhoods in~$\law\ust_t$ are at distance at least~$2s$ from one another.
If~$J\ust_u$ denotes the set of indices of the points of~$\tY$ that enter~$\law\ust$ up to time~$u$, 
then $\law_u\ust$ consists of disjoint neighbourhoods $(\KK(\bP_j,u-\bt_j),\, j\in J\ust_u)$, and
new points are generated at rate $\r \sum_{j\in J\ust_u}\vva_{u-\bt_j}(1-\p\ust_u)$, where the
censoring probability $\p\ust_u$ is given by
\eq\label{censoring-probability}
    \p\ust_u \Def L^{-1}\sum_{j\in J\ust_u}\vva_{2s+(t-\bt_j)+(t-u)}.
\en
In our applications, we can find suitably small bounds for~$\p\ust_u$, so that the growth of
the numbers of neighbourhoods in~$\law\ust$ is still reasonably close to that of the CMJ process~$\bX$.
In view of the `hard core' censoring, the points $(\bP_j,\,j\in J\ust_u)$ are no longer independent of
one another, but their marginal distribution is still uniform on~$\cC$ if~$P_0$ is chosen at random.
Note also that $\law\ust_u \subset \law_u$ for each $s,t \ge 0$ and $0 < u \le t$.

We shall also use comparisons
between the CMJ process~$\bX$ and `flattened' versions $\hX_-$, $\hX_0$ and~$\hX_+$ that
are of the form discussed in the previous section.  We start by noting that, from
the inequality~\Ref{AB-volume-approx}, 
\eq\label{AB-volume-bounds}
   \vva s^d\{1 - \etaL\} \Le \vva_s \Le \vva s^d\{1 + \etaL\},\quad 0 < s \le \tlm,
\en 
where~$\tlm := \frac3{2\l}\log\L$ is as in~\Ref{AB-tlm-def}, and
\eq\label{AB-eta-Lambda-def}
  \etaL \Def c_g\biggl( \frac{3\log\L}{2\L^{1/d}} \biggr)^{\g_g}.
\en
Hence, up to time~$\tlm$, the process~$\bX$ is stochastically dominated
by the flattened process~$\hX_+$, defined as in the previous section, having intensity 
$\r_+ := \r(1+\etaL)$ per unit volume, and hence growth rate
$\l_+ := \l\{1+\etaL\}^{1/d}$; similarly, it stochastically dominates the flattened process~$\hX_-$ 
with $\r_- := \r(1-\etaL)$ and $\l_- := \l\{1-\etaL\}^{1/d}$.  We also define the
flattened process~$\hX_0$ with intensity~$\r$ per unit volume, and with growth rate~$\l$.
The quantities $M_j^+$, $M_j^0$ and~$M_j^-$, and their standardized versions 
$H_j^+$, $H_j^0$ and~$H_j^-$, correspond to these processes.
We make the relationships between the processes precise with
the following construction.

\begin{lemma}\label{AB-dominating-processes}
Let the successive birth times in the branching processes $\bX$, $\hX_-$, $\hX_0$ and~$\hX_+$ be denoted
by $(\bt_j,\htau_j^-,\htau_j^0,\htau_j^+,\,j\ge0)$, respectively, and let $(T_t,T_t^-,T_t^0,T_t^+)$ denote the sets
of birth times up to time~$t$ in each of the processes.
If, for some $0\le s < \tlm$, $T_s^- \subset T_s \subset T_s^+$
and $T_s^- \subset T_s^0 \subset T_s^+$, then the processes
$\bX$, $\hX_-$, $\hX_0$ and~$\hX_+$ can be defined on the same probability space,
in such a way that, for all $s\le t \le \tlm$,
\[
    T_t^- \subset T_t \subset T_t^+ \quad \mbox{and}\quad T_t^- \subset T_t^0 \subset T_t^+ \quad\mbox{a.s.}
\]
\end{lemma}

\proof
The birth rate of~$\MBP$ at time~$t$ is given by
\[
    r(\MBP,t)  \Def \r \sum_{j\colon\t_j\in T_t} \vva_{t-\bt_j},
\]
and of $\hX_0$ by
\[
    r(\hX_0,t) \Def \l H_d^0(t) \Eq \l^{d+1} \sum_{j\colon\htau_j^0\in T_t^0} (t-\htau_j^0)^d/d!
     \Eq \r \vva \sum_{j\colon\htau_j^0\in T_t^0} (t-\htau_j^0)^d,
\]
with analogous representations for $r(\hX_-,t)$ and~$r(\hX_+,t)$.  Thus, for any time~$t$ such that
\eq\label{AB-inclusions}
   T_t^- \subset T_t \subset T_t^+ \quad \mbox{and}\quad T_t^- \subset T_t^0 \subset T_t^+,
\en
we have $r(\hX_-,t) \le r(\MBP,t) \le r(\hX_+,t)$ and $r(\hX_-,t) \le r(\hX_0,t) \le r(\hX_+,t)$.
Hence, for~$s$ as given, we can construct all four processes on the same probability space,
for $s \le t \le \tlm$, by realizing~$\hX_+$ on $[s,\tlm]$
together with an independent sequence of
independent random variables~$(U_j,\,j\ge1)$ uniformly distributed on $[0,1]$, and then
thinning in the following way.
At each successive point $\htau_j^+ > s$, include it as a point of~$\MBP$ if $U_j r(\hX_+,t) \le r(\MBP,t)$;
similarly, if $U_j r(\hX_+,t) \le r(\hX_-,t)$, include~$\htau_j^+$ as a point of~$\hX_-$, and
if $U_j r(\hX_+,t) \le r(\hX_0,t)$, include~$\htau_j^+$ as a point of~$\hX_0$.  This
construction preserves the inclusions~\Ref{AB-inclusions} for all times up to $\tlm$, 
and, because independently thinned Poisson processes are again Poisson processes, also
yields the right distributions for the processes $\MBP$, $\hX_0$ and~$\hX_-$.
\ep
  
\nin In what follows, we shall use~$\fpp t$ to denote the filtration for the combined
construction in Lemma~\ref{AB-dominating-processes}.  We shall henceforth only consider times
in $[0,\tlm]$, and will take~$\L$ large enough that 
\eq\label{Lambda-conds}
         \exp\{3\etaL\tlm\} \Le 2\quad\mbox{and}\quad \etaL \Le 1.
\en

\subsection{Relating the proportion informed to the function~$\ell$}
  The first step in our detailed calculations is to replace~$L_t/L$
with $\ex\{L_t/L \giv \tff_s\}$, where $\tff_s := \s(\tY_u,\,0 \le u\le s)$, for suitable $s<t$;  
this conditional expectation is easier to handle. We start by bounding the conditional variance 
$\var\{L_t/L \giv \tff_s\}$, for suitable values of~$s < t$.

The basis for our argument is given by the observations that
\eq\label{AB-basic-repn}
  \ex\{1 - L_t/L \giv \tff_s\} \Eq \pr[K \notin \law_t \giv \tff_s] \quad\mbox{and}\quad
  \ex\{(1 - L_t/L)^2 \giv \tff_s\} \Eq \pr[K,K' \notin \law_t \giv \tff_s],
\en
where $K$ and~$K'$ are chosen independently and uniformly in~$\cC$, implying that
\eq\label{AB-variance}
   \var\{L_t/L \giv \tff_s\} \Eq \pr[K,K' \notin \law_t \giv \tff_s] - \{\pr[K \notin \law_t \giv \tff_s]\}^2.
\en
On the other hand, 
\eq\label{AB-gossip-identity}
   \{K \notin \law_t\} \Eq \{\ggg_{t,s}^K \cap \law_s = \emptyset\},
\en
where~$\ggg_{t,s}^K$ denotes the set of all points at time~$s$ that, if informed, would inform~$K$
by time~$t$.  Now, for the gossip process, $\ggg_{t,s}^K$ is {\it independent\/} of~$\tff_s$, and
has the same distribution as~$\law_{t-s}$.  In view of~\Ref{AB-gossip-identity}, we  thus have
\eq\label{AB-gossip-prob}
  \pr[K \notin \law_t\giv \tff_s] \Eq \pr[\ggg_{t,s}^K \cap \law_s = \emptyset\giv \tff_s],
\en
where $\law_s$ is $\tff_s$-measurable and $\ggg_{t,s}^K$ is independent of~$\tff_s$,
and
\eq\label{AB-gossip-prob-2}
  \pr[K,K' \notin \law_t\giv \tff_s] \Eq \pr[\{\ggg_{t,s}^K \cap \law_s = \emptyset\} 
               \cap \{\ggg_{t,s}^{K'} \cap \law_s = \emptyset\} \giv \tff_s],
\en
with $\ggg_{t,s}^K$ and~$\ggg_{t,s}^{K'}$ independent of~$\tff_s$, but not of each other.
Indeed, in view of~\Ref{AB-variance}, it is the extent of their dependence that measures
$\var\{L_t/L \giv \tff_s\}$.

Writing $t_s := t-s$,  our argument now involves bounding the differences
\eqa\label{AB-prob-diff-1}
  &&\pr[\ggg_{t,s}^K \cap \law_s = \emptyset\giv \tff_s] - \pr[\bbb^K(t_s) \cap \law_s = \emptyset\giv \tff_s]
     \quad\mbox{and} \\
  &&\pr[\{\ggg_{t,s}^K \cap \law_s = \emptyset\} 
               \cap \{\ggg_{t,s}^{K'} \cap \law_s = \emptyset\} \giv \tff_s] \phantom{XXXXXXXXXXXXXXXXXXXXX}  \non\\ && \qquad
   \qquad\qquad\qquad\qquad\mbox{} - \pr[\{\bbb^K(t_s) \cap \law_s = \emptyset\} 
               \cap \{\bbb^{K'}\!(t_s) \cap \law_s = \emptyset\} \giv \tff_s]  \label{AB-prob-diff-2}
\ena 
between the probabilities \Ref{AB-gossip-prob}
and~\Ref{AB-gossip-prob-2} and the smaller ones obtained by replacing $\ggg_{t,s}^K$ and~$\ggg_{t,s}^{K'}$
by their related (independent) branching and growth processes $\baL^K$ and~$\baL^{K'}$.
\ignore{
derived from (independent) CMJ processes 
$\MBPK$ and~$\MBPKD$, each distributed as~$\bX$, by associating independent and uniformly distributed 
locations in~$\cC$ to their event times. Since their branching rates ignore overlap in~$\cC$,
}
These, as observed in the joint construction at the beginning of the section, give rise to
stochastically larger sets than $\ggg_{t,s}^K$ and~$\ggg_{t,s}^{K'}$.
If both of the differences \Ref{AB-prob-diff-1} and~\Ref{AB-prob-diff-2} are smaller than some~$\e$, 
then the independence of $\bbb^K$ and~$\bbb^{K'}$
immediately implies that $\var\{L_t/L \giv \tff_s\} \le 4\e$.  Using this strategy, we prove the following
lemma.

\begin{lemma}\label{AB-variance-lemma}
Under the above assumptions, there is a 
constant $C_{\ref{AB-variance-lemma}} = C_{\ref{AB-variance-lemma}}(d)$ such that
\[
     \var\{L_t/L \giv \tff_s\} \Le   C_{\ref{AB-variance-lemma}} 
                  \L^{-2}(1 + (\l s)^d) e^{2\l(t-s)}(\l^d M_s +  N_s).
\]
\end{lemma}

\proof
To control the differences \Ref{AB-prob-diff-1} and~\Ref{AB-prob-diff-2},   we begin by running a process~$\tY^K$, 
defined following~\Ref{ADB-Y-proc-def}, until time~$t_s$, and thin to obtain~$\ggg_{t,s}^K$. As in~\Ref{ADB-bJ-etc-defs}, 
let~$\bJ^K_u := \{j \ge 0\colon \bt_j^K \le u\}$, and set $\bN^K_u := |\bJ^K_u|$ and 
$\bM^K_u := \sum_{j\in \bJ^K_u}(u-\bt_j^K)^d$.  We then thin~$\tY^K$ further to construct the 
process~$(\lawtsK(u),\,0\le u\le t_s)$, by the method used to construct $\law\ust$ in~\Ref{ADB-lawst-rule}.

We now consider the difference 
\[
    \D_{s,t} \Def \pr[\lawtsK(t_s) \cap \law_s = \emptyset\giv \tff_s] 
                           - \pr[\bbb_{t_s}^K \cap \law_s = \emptyset\giv \tff_s],
\]
which is an upper bound for the real quantity~\Ref{AB-prob-diff-1} of interest to us.  The quantity~$\D_{s,t}$
is no larger than the conditional expectation given~$\tff_s$ of the number~$Z_{t,s}^K$ of intersections 
between {\it censored\/} islands of~$\baL^K_{t_s}$ and the islands of~$\law_s$.
If an island born in~$\bX^K$ at~$u$ is censored, the expected number of censored islands that result at~$t_s$
is at most $c_1e^{\l_+(t_s-u)}$, by~\Ref{AB-N-star-moments} and because~$\bX^K$ is stochastically
dominated by~$\hX_+$.  These islands each have radius at most $(t_s-u)$.
Hence, given~$\tff_s$, 
the expected number of 
intersections resulting from a censored island born at~$u$ is at most
\eqs
  \lefteqn{ c_1 e^{\l_+(t_s-u)} \sum_{j\in J_s}L^{-1} \vva_{(s-\t_j) + (t_s-u)}  }\\
        &&\Le  c_1 e^{\l_+(t_s-u)} \vva (1 + \etaL) L^{-1} \sum_{j\in J_s}((s-\t_j) + (t_s-u))^d \\
        &&\Le 2^d c_1 \vva  e^{\l_+(t_s-u)}L^{-1}(\bM_s + \bN_s(t_s-u)^d),
\ens
in view of \Ref{AB-intersection-prob}, \Ref{AB-volume-approx} and~\Ref{Lambda-conds};  
$\bN$ and~$\bM$ are as in~\Ref{ADB-bJ-etc-defs}.  
Similarly, using~\Ref{censoring-probability}, the conditional probability~$\p_u^{0,t_s,K}$ of an 
island born in~$\bX^K$ at~$u$ being censored for $\lawtsK$, given the history up to~$u$, 
is bounded above by
\eqs
    && (1+\etaL)\vva  L^{-1} \int_{(0,u)} \{2s + (t_s-v) + (t_s-u)\}^d \bN^K(dv) \\
        &&\qquad\Le 2.3^{d-1}\vva  L^{-1} \int_{(0,u)}  \{(2s)^d + (2(t_s-u))^d + (u-v)^d\} \bN^K(dv) \\[1ex]
        &&\qquad\Eq 2.3^{d-1}\vva  L^{-1} \bigl\{\bN^K_{u-}\{(2s)^d + (2(t_s-u))^d\} + \bM^K_{u-}\bigr\}.
\ens
Hence, again using~$\bN^K$ as an upper
bound for the number of uncensored islands, and noting that the birth intensity in~$\bX^K$ at time~$u$ is at most 
$$
          \r \sum_{j\in\bJ_u} \vva_{u-\bt_j^K} \Le 2\vva  \r \bM^K_u,
$$ 
we have
\eqa
  \lefteqn{ \ex\{Z_{t,s}^K \giv \tff_s\}} \non \\
  &\le&   \ex\Blb \int_0^{t_s} 2.3^{d-1}\vva L^{-1} \bigl\{\bN^K_{u-}
          \{ (2s)^d + 2^d(t_s-u)^d\} + \bM^K_{u-}\bigr\}\right.  \non\\
  &&\qquad\qquad\left. 2^d c_1\vva  e^{\l_+(t_s-u)} L^{-1}(\bM_s + \bN_s(t_s-u)^d)\,\bN^K(du) \Giv \tff_s \Brb \non \\
   &\le& 2^{d+1}3^{d-1} c_1 \r \{\vva \}^3  L^{-2} \label{AB-E-bnd-1}\\
    &&\ex\Blb \int_0^{t_s} \bigl\{\bN^K_u\{(2s)^d + 2^d(t_s-u)^d\} + \bM^K_u\bigr\} 
           e^{\l_+(t_s-u)}(\bM_s + \bN_s(t_s-u)^d)\bM^K_u\,du  \Giv \tff_s \Brb . \non
\ena
Now, by \Ref{AB-N-star-moments}, \Ref{AB-L-star-moments} and Cauchy--Schwarz, and because~$\bX^K$ is stochastically
dominated by~$\hX_+$,
$$
    \ex \bigl\{( \bN^K_u \{(2s)^d + (t_s-u)^d\} + \bM^K_u)\bM^K_u \bigr\} 
              \Le c_2 d! \l_+^{-d}\{ (2s)^d + (t_s-u)^d + d!\l_+^{-d}\} e^{2\l_+ u}.
$$  
Using this in~\Ref{AB-E-bnd-1}, and noting that $\l_+ \le \l(1+\etaL)$ and that $\r\vva d! = \l^{d+1}$,
gives the following bound for~\Ref{AB-prob-diff-1}:
\eqa
  0 &\le& \pr[\ggg_{t,s}^K \cap \law_s = \emptyset\giv \tff_s] - \pr[\bbb^K_{t_s} \cap \law_s = \emptyset\giv \tff_s] 
           \Le \ex\{Z_{t,s}^K \giv \tff_s\} \non \\
    &\le&   C_{1}(d)(1 + (\l_+s)^d) \l_+^{-d} \{\vva \}^2 L^{-2} e^{2\l_+ t_s}(\bM_s + \l_+^{-d}\bN_s) \\
    &\le&   C_{1}(d) \L^{-2} (1 + (\l_+s)^d) e^{2\l t_s}(\l^d \bM_s + \bN_s). \phantom{XXX}
    \label{AB-E-bnd-2}
\ena

We now need to bound~\Ref{AB-prob-diff-2}.
This can be done by introducing a process~$\law^{0,t_s,K,K'}$, constructed in the same way as~$\lawtsK$,
but starting from two initial points $K,K'$ and using a CMJ process~$\bX^{K,K'}$, which is the same as
using two independent CMJ processes $\bX^K$ and~$\bX^{K'}$, by the branching property.  
Now $\law^{0,t_s,K,K'}(t_s) \subset (\ggg_{t,s}^K\cup\ggg_{t,s}^{K'})$, and the conditional expection 
given~$\tff_s$ of the number~$Z_{t,s}^{K,K'}$ of intersections between censored islands 
of~$\bX^{K,K'}_{t_s}$ and the islands of~$\law_s$ satisfies
\eq\label{AB-E-bnd-3}
   \ex\{Z_{t,s}^{K,K'} \giv \tff_s\} \Le C_2(d) \{\vva \}^2 L^{-2} (1 + (\l_+s)^d) e^{2\l_+ t_s}(\l^d\bM_s + \bN_s),
\en
by an argument exactly as before, but
for a larger constant~$C_2(d)$ than~$C_1(d)$ appearing in~\Ref{AB-E-bnd-2}.  
Since~$\ex\{Z_{t,s}^{K,K'} \giv \tff_s\}$ 
is a bound for the difference in~\Ref{AB-prob-diff-2}, we have enough to prove the lemma.
\ep

\begin{remark}
With $s = \a_1\l^{-1}\log \L$ and $t = \a_2\l^{-1}\log \L$,
where $\a_1 < \a_2 \le 1$, and since $\ex(\l^d\bM_s+\bN_s) =  O(e^{\l_+ s})$, from~\Ref{AB-L-star-moments},
it follows that $\var\{L_t/L \giv \tff_s\}$ is typically of order $O\bigl(\L^{2\a_2-\a_1-2}(\log\L)^d \bigr)$.
\end{remark}

\medskip
Our main interest is in approximating the distribution of~$L_t/L$ when 
\eq\label{AB-t-Lambda-def}
  t \Eq t_\L(u) \Def \l^{-1}\{\log\L + u\},
\en 
for~$u$ fixed.  This is because the times $(t_\L(u),\,u\in\re)$ asymptotically represent 
the period in which $L_t/L$ increases from $0$ to~$1$. 
Taking $\a_1 = \a < 1$ and $\a_2 = 1$ in  the remark, it follows that
$\var\{L_{t_\L(u)}/L \giv \tff_s\}$ is typically of
order $O(\L^{-\a})$ for $s := \a\l^{-1}\log \L$.  Now pick $v := \a_1\l^{-1}\log \L$  and $s := \a_2\l^{-1}\log\L$,
with $\a_1 < \a_2 < 1$.  Then
\[
     \var\{L_t/L \giv \tff_v\} \Eq \var\{\ex(L_t/L \giv \tff_s) \giv \tff_v\}
          + \ex\{\var(L_t/L \giv \tff_s) \giv \tff_v\} \,,
\]
in which the latter term, again by the remark, is typically of order $O(\L^{-\a_2})$ if $t = t_\L(u)$.  
Supposing
that $\var\{L_t/L \giv \tff_v\}$ is actually of magnitude $\L^{-\a_1}$,  this indicates that the
conditional distribution of~$L_t/L$ given~$\tff_v$ is essentially that of the conditional
distribution of $\ex(L_t/L \giv \tff_s)$ given~$\tff_v$.  So the next step is
to examine  $\ex\{(1 - L_t/L) \giv \tff_s\}$ in detail, for $t = t_\L(u)$, and to 
express it in more amenable form.

The next lemma once again uses the backward branching process~$\baL^K$ from a randomly chosen point~$K$.
We define
$\ff^K_{s,t} := \tff_s \bigvee \ff^K_{t-s;0}$, where $\ff^K_{v;0} := \s(\bN^K_u,\,0\le u\le v)$
contains the information about when the islands of~$\baL^K$ were formed, up to time~$v$, but not where 
they are centred.  We then write~$Z^{s,t}$ for the number of islands of~$\baL^K_{t_s}$ that intersect~$\law_s$.

\begin{lemma}\label{AB-first-representation}
With the definitions above, there is a 
constant $C_{\ref{AB-first-representation}} = C_{\ref{AB-first-representation}}(d)$ such that
\eqs
 \lefteqn{\bigl|\ex\{(1 - L_t/L) \giv \tff_s\} - \ex\bigl\{\exp\{-\mmm^K_{s,t}\} \giv \tff_s \bigr\}\bigr| }\\
    &&\Le  C_{\ref{AB-first-representation}}
        \{\L^{-1} \bN_s(\l t)^d   +  \L^{-2}(1 + (\l_+s)^d)e^{2\l(t-s)}(\l^d \bM_s + \bN_s)\},
\ens
where $\mmm^K_{s,t} := \ex\{Z^{s,t} \giv \ff^K_{s,t}\}$.
\end{lemma}

\proof
We start by using \Ref{AB-basic-repn}, \Ref{AB-gossip-identity}
and~\Ref{AB-E-bnd-2} to show that, for $t>s$,
\eq\label{AB-simplify-mean-1}
   | \ex\{(1 - L_t/L) \giv \tff_s\} -  \pr[\bbb^K_{t_s} \cap \law_s = \emptyset\giv \tff_s]|
        \Le C_{1}(d) \L^{-2} (1 + (\l_+s)^d) e^{2\l t_s}(\l^d \bM_s + \bN_s).
\en
We now use Poisson approximation to approximate the probability 
$\pr[\bbb^K_{t_s} \cap \law_s = \emptyset\giv \tff_s]$,
using the conditional independence between the locations of
the islands of~$\bbb^K_{t_s}$, given~$\ff^K_{s,t}$, as the basis of the approximation.  


We first observe that the conditional probability that an island of~$\bbb^K_{t_s}$ with radius~$v$ intersects~$\law_s$, 
given~$\ff^K_{s,t}$, is at most
\eq\label{AB-intersect-1}
     \sum_{j\in J_s} \vva_{s-\t_j+v}L^{-1}   \Le 2 \bN_s \vva L^{-1}t^d 
            \Eq (1+\etaL)\L^{-1}\bN_s(\l t)^d ,
\en
in view of~\Ref{AB-intersection-prob}, by \Ref{AB-volume-approx}, \Ref{AB-eta-Lambda-def} and~\Ref{Lambda-conds},
and because $v \le t-s$.
This, using~$Z^{s,t}$ to denote the number of islands of~$\baL^K_{t_s}$ that intersect~$\law_s$, implies that
\eq\label{AB-Poisson-1}
    \dtv(\law(Z^{s,t} \giv \ff^K_{s,t}),\Po(\mmm^K_{s,t})) \Le 2 \L^{-1}\bN_s(\l t)^d,
\en  
by \BHJ~(1992, (1.23)), where $\mmm^K_{s,t} := \ex\{Z^{s,t} \giv \ff^K_{s,t}\}$.
Hence, from~\Ref{AB-Poisson-1},
\[
    \bigl|\pr[Z^{s,t} = 0 \giv \tff_s] - \ex\{\exp(-\mmm^K_{s,t}) \giv \tff_s\}\bigr|
       \Le 2\L^{-1}\bN_s(\l t)^d,
\]
and combining this with~\Ref{AB-simplify-mean-1} gives the lemma. 
\ep

We now define
\eq\label{AB-tilde-M-def}
    \mmc^K_{s,t} \Def \int_0^{t_s} \sum_{j \in J_s} \vva  L^{-1}(s - \t_j + t_s-v)^d \, \bN^K(dv),
\en
as an approximation to~$\mmm^K_{s,t}$.
The following lemma bounds the accuracy of the approximation
for $t = t_\L(u)$.

\begin{lemma}\label{AB-M-approx}
For any $\gda > 0$, there is an event $B_{\ref{AB-M-approx}}(\gda,s) \in \tff_s$ with 
$\pr[\{B_{\ref{AB-M-approx}}(\gda,s)\}^c] \le C_{\ref{AB-M-approx}}\L^{-\gda}$  
such that, for $t = t_\L(u)$,
\[
   \ex\{|\mmc^K_{s,t} - \mmm^K_{s,t}|\giv \tff_s\} I[ B_{\ref{AB-M-approx}}(\gda,s)] 
             \Le C'_{\ref{AB-M-approx}} \L^{\gda} e^u \{\L^{-1} \{\l s\}^d e^{\l s} + \etaL\}, 
\]
where $C_{\ref{AB-M-approx}}$ and~$C'_{\ref{AB-M-approx}}$ depend only on~$d$.
\end{lemma}

\proof
We begin by introducing the censored version $\bbb\uss$ of the process~$\bbb$.
We denote the indices of islands in~$\bbb\uss_s$ by $\Js \subset J_s$, and write $r_{js} := s - \bt_j$.
It then follows that
\eq\label{AB-MM-inequalities}
   \int_0^{t_s} \sum_{j\in \Js} L^{-1} \vva_{r_{js}+t_s-v} \bN^K(dv)
         \Le \mmm^K_{s,t} \Le  \int_0^{t_s}\sum_{j\in \bJ_s} L^{-1} \vva_{r_{js}+t_s-v} \bN^K(dv),
\en
with the lower bound using the separation between the islands of~$\bbb\uss$.
Now, from \Ref{AB-volume-bounds}, \Ref{AB-eta-Lambda-def} and~\Ref{AB-MM-inequalities},
\eqs
  \mmm^K_{s,t} &\ge& \int_0^{t_s} \sum_{j \in \Js} L^{-1} \vva_{r_{js}+t_s-v} \bN^K(dv) \\  
      &\ge& (1-\etaL)\int_0^{t_s} \sum_{j \in \Js} L^{-1} \vva (r_{js}+t_s-v)^d \bN^K(dv),
\ens
and 
\eqs
   \mmm^K_{s,t} &\le& \int_0^{t_s}\sum_{j\in \bJ_s} L^{-1} \vva_{r_{js}+t_s-v} \bN^K(dv) \\
      &\le& (1+\etaL) \int_0^{t_s} \sum_{j \in \bJ_s} L^{-1} \vva (r_{js}+t_s-v)^d \bN^K(dv).
\ens
Hence 
\eq\label{ADB-new-bigger-piece}
    \mmm^K_{s,t} - \mmc^K_{s,t} \Le  \etaL \mmc^K_{s,t}
       + (1+\etaL)\int_0^{t_s} \sum_{j \in \bJ_s\setminus J_s} L^{-1}\vva (r_{js} + t_s-v)^d \, \bN^K(dv),
\en
and
\eq\label{AB-bigger-piece}
    \mmc^K_{s,t} - \mmm^K_{s,t} \Le \etaL \mmc^K_{s,t} 
           + (1+\etaL)\int_0^{t_s} \sum_{j \in J_s\setminus \Js} L^{-1}\vva (r_{js} + t_s-v)^d \, \bN^K(dv).
\en
This implies that
\eqa
   |\mmc^K_{s,t} - \mmm^K_{s,t}| &\le& \etaL \mmc^K_{s,t}
     + (1+\etaL)\int_0^{t_s} \sum_{j \in \bJ_s\setminus \Js} L^{-1}\vva (r_{js} + t_s-v)^d \, \bN^K(dv)\non\\
        &\le& \etaL \mmc^K_{s,t} + 
           2^{d}  \Bl \bN^K_{t_s}\sum_{j \in \bJ_s\setminus \Js} L^{-1} r_{js}^d +
                    L^{-1}(\bN_s - \nnn\uss_s) \bM^K_{t_s} \Br ,\label{AB-Poisson-mean-2}
\ena
where  $\nnn\uss_s := |\Js|$.
Thus we need to bound the conditional expectation given~$\tff_s$ of the right hand side of~\Ref{AB-Poisson-mean-2}.

Define $B_1(\gda,s)$ by
\eq\label{AB-B-def}
  B_1(\gda,s) \Def \Blb \l^d\sum_{j \in \bJ_s\setminus \Js} r_{js}^d + d!(\bN_s - \nnn\uss_s) 
             \Le \L^{-1+\gda} \{\l_+ s\}^d e^{2\l_+ s}\Brb \ \in\ \tff_s.
\en
Since $\baL^K$ is independent of~$\law$ in~\Ref{AB-Poisson-mean-2}, it follows that we can easily take the expectation,
given~$\tff_s$, of its second term.
For $t = t_\L(u)$, and using \Ref{AB-N-star-moments} and~\Ref{AB-L-star-moments}, this gives
\eqa
   \lefteqn{ \ex\Bigl\{ 2^{d}  \Bl \bN^K_{t_s}\sum_{j \in \bJ_s\setminus \Js} L^{-1} r_{js}^d +
                    L^{-1}(\bN_s - \nnn\uss_s) \bM^K_{t_s} \Br \Giv \tff_s \Bigr\} I[B_1(\gda,s)]} \non\\
 &&\Le
   2^d c_1 \L^{-1} e^{\l_+ t_s}\Bl \l^d\sum_{j \in \bJ_s\setminus \Js} r_{js}^d +
                    d!(\bN_s - \nnn\uss_s)  \Br I[B_1(\gda,s)]\phantom{XX} \non\\
 &&\Le 2^{d+2} c_1 \L^{-1+\gda} e^u \{\l_+ s\}^d e^{\l s}, \label{AB-Poisson-mean-diff}
\ena
where we have twice used $e^{(\l_+-\l)t} \le 2$ for $t \le \tml$, as follows from~\Ref{Lambda-conds}.
For the first term in~\Ref{AB-Poisson-mean-2}, from~\Ref{AB-tilde-M-def}, we have
\[
    \mmc^K_{s,t} \Le 2^d\vva L^{-1}\Bl \bN^K_{t_s}\sum_{j \in \bJ_s} r_{js}^d + \bN_s \bM^K_{t_s} \Br.
\]
Defining
\[
    B_2(\gda,s) \Def \Blb \l^d\sum_{j \in \bJ_s} r_{js}^d + d!\bN_s
             \Le \L^{\gda} e^{\l_+ s}\Brb,
\] 
it thus follows from the independence of $\law$ and~$\baL^K$ that, for $t = t_\L(u)$,
\eq\label{AB-Poisson-mean-diff-2}
   \etaL \ex\{ \mmc^K_{s,t} \giv \tff_s \}I[ B_2(\gda,s)] \Le 2^{d+1} d! c_1 \etaL \L^{\gda}e^u, 
\en
using~\Ref{Lambda-conds} to bound $e^{(\l_+ - \l)t}$.

To complete the proof of the lemma, we need to show that 
$$
       \pr[(B_1(\gda,s))^c] + \pr[(B_2(\gda,s))^c] \Eq O(\L^{-\gda}).
$$
For $\pr[(B_1(\gda,s))^c]$, we bound  $\ex\{\bN_s - \nnn\uss_s\}$ and 
$\ex\bigl\{ \sum_{j \in \bJ_s\setminus \Js}  r_{js}^d \bigr\}$, and then use Markov's inequality.
We begin by bounding the conditional probability~$\p\uss_u$, given the past up to time~$u- < s$, that an island of~$\bX$, 
born to an uncensored parent at~$u$, is censored in~$\bbb\uss$. Using~\Ref{censoring-probability}, it is no greater than
\[
    L^{-1}\int_{(0,u)}\vva_{s-v + s-u + 2s}^d\,\bN(dv) \Le (1+\etaL) \vva  L^{-1}(4s)^d \bN_{u-} \,.
\]
If it is censored, bounding~$\bX^K$ by the branching process~$\hX_+$ and  using \Ref{AB-N-star-moments} 
and~\Ref{AB-L-star-moments}, the expected number 
of its offspring by time~$s$, all of which are also censored, is at most
$c_1e^{\l_+(s-u)}$, and the expected volume censored at most $c_1 d! \l_+^{-d} e^{\l_+(s-u)}$.
Hence
\eqa
    \lefteqn{\ex\{\bN_s - \nnn\uss_s\}} \non\\
        &&\Le (1+\etaL) \vva  L^{-1}(4s)^d \,
                  \ex \Blb \int_0^s c_1 e^{\l_+(s-u)} \bN_{u-}\,\bN(du) \Brb \non \\
        &&\Le c_1 (1+\etaL) \vva  L^{-1}(4s)^d \,
                  \ex \Blb \int_0^s e^{\l_+(s-u)} M_0^+(u) \r(1+\etaL) + \vva  M_d^+(u)\,du \Brb \non\\
        &&\Le 4c_1c_2 \r \{\vva \}^2 d!\l_+^{-d} L^{-1}(4s)^d \int_0^s e^{\l_+(s+u)}\,du \non\\
        &&\Le 4c_1c_2 (1+\etaL) \L^{-1}(4\l s)^d e^{2\l_+ s}, \label{AB-N-diff-bnd}
\ena
again by \Ref{AB-N-star-moments} and~\Ref{AB-L-star-moments}, and from Cauchy--Schwarz.
Then, by a similar argument,
\eqa
    \ex\Blb \sum_{j \in \bJ_s\setminus \Js}  r_{js}^d \Brb  
      &\le& (1+\etaL) \vva  L^{-1} (4s)^d \, \ex \Blb \int_0^s c_1 d!\l^{-d} e^{\l_+(s-u)} \bN_{u-}\,\bN(du) \Brb \non\\
      &\le&  2c_1c_2 d! \L^{-1} (4s)^d e^{2\l_+ s}. \label{AB-N-sum-bnd}
\ena
Combining \Ref{AB-N-diff-bnd} and~\Ref{AB-N-sum-bnd} 
and using Markov's inequality, $\pr[\{B_1(\gda,s)\}^c] \le c\L^{-\gda}$,
for a constant~$c$ depending only on~$d$.

For $\pr[(B_2(\gda,s))^c]$, we again bound~$\bX^K$ by the branching process~$\hX_+$ and use \Ref{AB-N-star-moments} 
and~\Ref{AB-L-star-moments}, giving
\eq\label{AB-N-growth}
   \ex \bN_s \Le c_1 e^{\l_+ s};\qquad
    \ex\Blb \sum_{j \in \bJ_s}  r_{js}^d \Brb \Le \ex M^+_d(s) \Le c_1 d! \l_+^{-d}e^{\l_+ s}.
\en
Hence, from Markov's inequality, $\pr[\{B_1(\gda,s)\}^c] \le c'\L^{-\gda}$,
for a constant~$c'$ depending only on~$d$, and the lemma is proved by taking $B_0(\gda,s) = B_1(\gda,s) \cap B_2(\gda,s)$.
\ep

We now replace $\ex\{\exp(-\mmc^K_{s,t}) \giv \tff_s\}$ by an expression involving the function~$\ell$ defined
in~\Ref{ADB-l-def}, and using the quantity~$W^*(s)$ defined by
\eq\label{ADB-W-star-def}
  W^*(s) \Def e^{-\l s}\sum_{l=0}^d \sum_{j\in \bJ_s} \frac{(\l (s-\bt_j))^l}{l!}
         \Le e^{-\l s} \sum_{l=0}^d H_l^+(s) \Eq e^{(\l_+-\l)s} W^+(s),  
\en
where the inequality follows from
Lemma~\ref{AB-dominating-processes}, so that, from~\Ref{Lambda-conds},
$\ex W^*(s) \Le 2$.

\begin{lemma}\label{AB-e-M-approx}
Take $s \le \l^{-1}\log\L$, and
let $\mmc^K_{s,t}$ be defined as in~\Ref{AB-tilde-M-def}, $W^*(s)$ as in~\Ref{ADB-W-star-def} and~$\ell$ as
for Lemma~\ref{ADB-l-def}.  Then, for any $\gda > 0$ and $0 < \h < \rrr(d)$, there is an 
event~$B_{\ref{AB-e-M-approx}}(\gda,\h,s) \in \tff_s$ and constants $C_{\ref{AB-e-M-approx}}$ and~$C'_{\ref{AB-e-M-approx}}$,
depending only on~$d$, such that 
$$
   \pr[\{B_{\ref{AB-e-M-approx}}(\gda,\h,s)\}^c] \Le C_{\ref{AB-e-M-approx}}(\L^{-\gda} + \l s e^{-2\l (\rrr(d)-\h) s})
$$ 
and that
\eqs
 \lefteqn{\bigl|\ex\{e^{-\mmc^K_{s,t_\L(u)}} \giv \tff_s\}  - \bigl(1 - \ell(\log[\hc_d W^*(s)] + u)\bigr)\bigr|\,
            I[B_{\ref{AB-e-M-approx}}(\gda,\h,s)] }\\
     &&\Le C'_{\ref{AB-e-M-approx}}(1+e^u)\bigl(\L^{\gda}(\etaL\log \L + \L^{-1}e^{\l s}) + e^{-\l\h s}\bigr),
\ens
uniformly in $t_\L(u) \le \tml$,  where $\hc_d := d!/(d+1)$. 
\end{lemma}
 
\proof
We first observe, from \Ref{AB-tilde-M-def}  and~\Ref{AB-Lambda-def} that 
\eqa
   \mmc^K_{s,t} &=& L^{-1}\vva \sum_{j\in J_s} 
             \int_0^{t_s} \sum_{l=0}^d {d \choose l} r_{js}^l (t_s-u)^{d-l} \bN^K(du)  \non\\
   &=& \L^{-1}  \sum_{l=0}^d {d \choose l} \Bl \sum_{j\in J_s} \{\l r_{js}\}^l \Br
      \,\Bl \int_0^{t_s} \{\l(t_s-u)\}^{d-l} \bN^K(du) \Br,  \label{AB-M-rewrite}
\ena
with $r_{js} := s - \bt_j$ as before.
Now realize $\hX_-$, $\MBP$ and~$\hX^+$ together as in Lemma~\ref{AB-dominating-processes}, so that 
\eq\label{AB-sandwich}
    H_l^-(s) \Le  \sum_{j\in J_s} \frac{(\l r_{js})^l}{l!} \Le H_l^+(s)\ \mbox{a.s.,\quad for}\ 
                 0\le s\le \tlm.
\en
Then, for such~$s$, it follows from \Ref{AB-H-sum-approx}, then using Lemma~\ref{ADB-phi-smooth}, \Ref{AB-E22-def} 
and~\Ref{AB-E22-bnd}, that, on an event 
$B_1^+(\h,s) \in \fpp s$ such that $\pr[\{B_1^+(\h,s)\}^c] \le c(d)(1+\l s) e^{-2\l (\rrr(d)-\h)s}$, 
we have
\eq\label{AB-C-upper-bound}
    \sum_{l=0}^d c_l H_l^+(s) \Le  C(1)\Bl \frac{1}{d+1}\sum_{l=0}^d H_l^+(s) + e^{\l_+(1-\h) s} \Br
\en
and
\eq\label{AB-C-lower-bound}
    \sum_{l=0}^d c_l  H_l^-(s) \ \ge\  C(1)\Bl \frac{1}{d+1}\sum_{l=0}^d H_l^-(s) - e^{\l_-(1-\h) s} \Br,
\en
for {\it all\/} choices of $c_0,\ldots,c_d$, where $C(1) := \sum_{l=0}^d c_l$.
Define 
\eq\label{AB-B3-def}
   B_2^+(\gda,s) \Def \Blb \frac{1}{d+1}\sum_{l=0}^d (H_l^+(s) - H_l^-(s)) \le e^{\l s} \L^{\gda} \etaL \log\L \Brb
    \ \in\ \fpp s.
\en
Then, on $B_1^+(\h,s) \cap B_2^+(\gda,s)$ and for $0\le s\le \tlm$, 
we have
\eqa
  \sum_{j\in J_s} \frac{(\l r_{js})^l}{l!} &\le& H_l^+(s) \Le \frac{1}{d+1}\sum_{l=0}^d H_l^+(s) + e^{\l_+(1-\h) s}
               \non\\
  &\le& \frac{1}{d+1}\sum_{l=0}^d H_l^-(s) + e^{\l_+(1-\h) s} + e^{\l s} \L^{\gda} \etaL \log\L \non\\
  &\le& \frac{1}{d+1} e^{\l s} W^*(s) + \e(\gda,\h,s), \label{AB-B1B2-consequence}
\ena
for all $0 \le l\le d$, 
from \Ref{AB-sandwich}, \Ref{AB-C-upper-bound}, \Ref{ADB-W-star-def} and~\Ref{AB-B3-def}, where 
\[
   \e(\gda,\h,s) \Def 2e^{\l(1-\h) s} + e^{\l s} \L^{\gda} \etaL \log\L.
\]
Arguing analogously, we also deduce that 
\[
   \sum_{j\in J_s} \frac{(\l r_{js})^l}{l!} \ \ge\  \frac{1}{d+1} e^{\l s} W^*(s) - \e(\gda,\h,s).
\]

Now $\pr[\{B_1^+(\h,s)\}^c] \le c(d)(1+\l s) e^{-2\l (\rrr(d)-\h)s}$.  Then, since 
$$
          \ex\Blb \sum_{l=0}^d H_l(s) \Brb \Eq e^{\l s}\ex W(s) \Eq e^{\l s}, 
$$
and using~\Ref{Lambda-conds}, we have
\[
   \ex\Blb \sum_{l=0}^d (H_l^+(s) - H_l^-(s)) \Brb \Eq e^{\l_+ s} - e^{\l_- s} \Le 8e^{\l s}\etaL \log\L
\]
in $0\le s\le \tlm$, and hence, by Markov's inequality,
\eq\label{AB-b2-prob}
   \pr[\{B_2^+(\gda,s)\}^c] \Le 8\L^{-\gda}.
\en   
Thus the event
\eq\label{AB-B-star-def}
  B_3(\gda,\h,s) \Def 
    \bigcap_{l=0}^d \Blb \Blm \sum_{j\in J_s} \frac{(\l r_{js})^l}{l!} - \frac{1}{d+1} e^{\l s} W^*(s) \Brm
     \le \e(\gda,\h,s) \Brb \ \in\ \tff_s
\en
is such that
\eq\label{B3c-prob}
   \pr[\{B_3(\gda,\h,s)\}^c] \Le C_1(d)(\L^{-\gda} + \l s e^{-2\l (\rrr(d)-\h)s}),
\en
for a suitable constant~$C_1(d)$.

Now, taking $c_l := \L^{-1}C_l(s,t)$, where
\eq\label{AB-c-choice}
  C_l(s,t) \Def  \int_0^{t_s} \frac{d!\{\l(t_s-u)\}^{d-l}}{(d-l)!} \bN^K(du),
\en
\Ref{AB-M-rewrite} implies that
\eqa
   \lefteqn{ \Blm \mmc^K_{s,t} -  \frac{e^{\l s}W^*(s)}{\L(d+1)} \sum_{l=0}^d C_l(s,t) \Brm I[B_3(\gda,\h,s)] }\non\\
   && \qquad\Le \L^{-1}\sum_{l=0}^d C_l(s,t)  e^{\l s} \{2e^{-\l\h s} + \L^{\gda}\etaL \log\L  \}.\non 
\ena
Hence also
\eqa
 \lefteqn{ \Blm \ex\Bl e^{-\mmc^K_{s,t}} \giv \tff_s \Br 
          - \ex\Bl \exp\Blb -  \frac{e^{\l s}W^*(s)}{\L(d+1)} \sum_{l=0}^d C_l(s,t)\Brb \Giv \tff_s\Br \Brm 
            I[B_3(\gda,\h,s)] } \non\\ 
    &&  \qquad  \Le\L^{-1}\ex\Blb\sum_{l=0}^d C_l(s,t) \Giv \tff_s \Brb  e^{\l s} \{2e^{-\l\h s} + \L^{\gda}\etaL \log\L  \}. 
           \phantom{XXXXXXXXX}
     \label{AB-Mst-approx-2'}
\ena

Now, because~$\bX^K$ can also be bounded between copies $\hX_-^K$ and~$\hX_+^K$ of $\hX_-$ and~$\hX_+$, 
using Lemma~\ref{AB-dominating-processes}, we have the inequality
\eq\label{AB-HK-sandwich}
    d!  H_{d-l}^{K,-}(t_s) \Le C_l(s,t) \Le d! H_{d-l}^{K,+}(t_s), \qquad 0 \le l \le d.
\en
Hence, since the $K$-processes can be chosen to be independent of~$\tff_s$, it follows that
\eq\label{AB-C-means}
  \ex\Blb \sum_{l=0}^d C_l(s,t)  e^{\l s} \giv \tff_s \Brb \Le 
       d! e^{\l_+(t-s) + \l s} \ex\{W^K(t_s)\} \Le 2\,d! e^{\l t},
\en
for any $0 < s \le t \le \tlm$.  Thus, from~\Ref{AB-Mst-approx-2'}, it follows that
\eqa
 \lefteqn{ \Blm \ex\Bl e^{-\mmc^K_{s,t}} \giv \tff_s \Br 
          - \ex\Bl \exp\Blb -  \frac{e^{\l s}W^*(s)}{\L(d+1)} \sum_{l=0}^d C_l(s,t)\Brb \Giv \tff_s\Br \Brm 
            I[B_3(\gda,\h,s)] } \non\\ 
    &&  \qquad  \Le  2\,d! \L^{-1}  e^{\l t} \{2e^{-\l\h s} + \L^{\gda}\etaL \log\L  \}. 
           \phantom{XXXXXXXXXXXXXXXXXXXX}
     \label{AB-Mst-approx-2}
\ena

The next step is to examine the difference 
$$
   \Blm \ex\Bl \exp\Blb -  \frac{e^{\l s}W^*(s)}{\L(d+1)} \sum_{l=0}^d C_l(s,t)\Brb \Giv \tff_s\Br
          - \f^1_{\l(t_\L(u)-s)}(\hc_d e^u W^*(s))\Brm,
$$  
where $\f^1_s(\th) := \ex\{e^{-\th W^1(s)}\}$.  To start with,
from~\Ref{AB-HK-sandwich} and Lemma~\ref{ADB-martingales}, 
\[
    d! e^{\l_- t_s}W^{K,-}(t_s) \Le   \sum_{l=0}^d C_l(s,t) \Le d! e^{\l_+ t_s}W^{K,+}(t_s).
\]
Hence, for any non-negative and $\tff_s$-measurable random variable~$\Th_s$, we have
\eq\label{AB-K-sandwich}
   \f^+_{t_s}(\Th_s d! e^{\l_+ t_s})
            \Le \ex\Blb \exp\Bl-\Th_s \sum_{l=0}^d C_l(s,t)\Br \Giv \tff_s \Brb \Le
  \f^-_{t_s}(\Th_s d! e^{\l_- t_s}),
\en
where
\eq\label{AB-f-defs}
   \f^+_t(\th) \Def \ex\{e^{-\th W^+(t)}\} \Eq \f^1_{\l_+ t}(\th)\quad\mbox{and}\quad
   \f^-_t(\th) \Def \ex\{e^{-\th W^-(t)}\} \Eq \f^1_{\l_- t}(\th),
\en
and $\f^1$ is as above, with the final equalities a consequence of~\Ref{AB-standard-W}.
Since $\l(1-\etaL) \le \l_- \le \l_+ \le \l(1+\etaL)$, we conclude from Lemma~\ref{ADB-phi-smooth} 
and~\Ref{Lambda-conds} that
\eqa
  \max\{|\f^+_t(\th e^{\l_+ t}) - \f^+_t(\th e^{\l t})|, |\f^-_t(\th e^{\l_- t}) - \f^-_t(\th e^{\l t})|\}
           &\le& 2e^{-1}  \etaL \l t ;\non\\
  \max\{|\f^+_t(\th) - \f^1_{\l t}(\th)|,|\f^-_t(\th) - \f^1_{\l t}(\th)|\} &\le& \th e^{-1} \etaL\,\l t e^{-\l t},
     \label{AB-phi-diffs} 
\ena
as long as $t \le \tml$.
Taking $\Th(s) := \{(d+1)\L\}^{-1}e^{\l s}W^*(s)$ and $t = t_\L(u)$, and using \Ref{AB-K-sandwich}, 
\Ref{AB-phi-diffs} and~\Ref{Lambda-conds}, this gives
\eqa
\lefteqn{
    \Blm \ex\Blb\exp\Bl - \frac{e^{\l s}W^*(s)}{\L(d+1)} \sum_{l=0}^d C_l(s,t)\Br \Giv \tff_s \Brb
        -  \f^1_{\l(t-s)}( \hc_d e^u W^*(s)) \Brm} \non\\[1ex]
   &&\Le 4e^{-1} \etaL \l t_s +  e^{-1}\Th(s) d! e^{\l_+ t_s} \etaL\,\l t_s e^{-\l t_s}
               \phantom{XXXXXXXXX}\non\\[1ex]
   &&\Le 4e^{-1} \etaL (\log\L + u) +   3\hc_d e^{\l s} W^*(s) \L^{-1}\etaL \log\L.
          \label{AB-f1-approx}
\ena
From \Ref{ADB-W-star-def}, we have $\ex\{W^*(s)\} \le 2$.  Thus, defining
\eq\label{AB-B4-def}
   B_4(\gda,s) \Def \Blb W^*(s) \le \L^{\gda} \Brb \ \in\ \tff_s,
\en
it follows that $\pr[\{B_4(\gda,s)\}^c] \le 2\L^{-\gda})$, and, combining \Ref{AB-Mst-approx-2} 
and~\Ref{AB-f1-approx}, that
\eqa
   \lefteqn{ \Bigl|\ex\{e^{-\mmc^K_{s,t_\L(u)}} \giv \tff_s\}  
            - \f^1_{\l(t_\L(u)-s)}( \hc_d e^u W^*(s))\Bigr| I[B_3(\gda,\h,s) \cap B_4(\gda,s)] } \non\\
              &&   \Le  C_2(d)(\L^{\gda}\etaL\log \L +  e^{-\l\h s}\L^{-1}e^{\l s}),\phantom{XXXXXXXXXXXXXXXX}
                       \label{Final-split}
\ena
uniformly in $t_\L(u) \le \tml$.  But now, from Lemma~\ref{ADB-phi-smooth}, on the event~$B_4(\gda,s)$,
\eqs
    |\f^1_{\l(t_\L(u)-s)}( \hc_d e^u W^*(s)) - \f^1_\infty( \hc_d e^u W^*(s))| 
             &\le& \frac1{2e}\hc_d e^u W^*(s) \exp\{-\l(t_\L(u)-s)\} \\
             &\le& \frac1{2e}\hc_d \L^{\gda-1}e^{\l s}, 
\ens
and $\f^1_\infty( \hc_d e^u W^*(s)) = 1 - \ell(\log( \hc_d W^*(s)) + u)$ by~\Ref{ADB-l-def},
\Ref{AB-standard-W} and~\Ref{W-convergence}.
This establishes the lemma, with
 $B_{\ref{AB-e-M-approx}}(\gda,\h,s) := B_3(\gda,\h,s) \cap B_4(\gda,s)$,
in view of \Ref{B3c-prob} and~\Ref{AB-B4-def}.
\ep

\subsection{Replacing $W^*(s)$ by $W(s,v)$}
Our aim is to approximate the conditional distribution of~$L_{t_\L(u)}/L$, given~$\tff_v$,
for suitably chosen~$v$.  After Lemma~\ref{AB-e-M-approx}, the problem has largely been reduced to
considering the conditional distribution of~$W^*(s)$.  However, in order to use the results of
Section~\ref{branching}, it is advantageous to replace $W^*(s)$ by a function of a {\it flattened\/}
branching process;  $W^*(s)$ is constructed from the birth times~$\bt_j$ of the original branching
process~$\bX$.  Accordingly, we define 
\eq\label{AB-Wsv-def}
     W(s,v) \Def e^{-\l s} \sum_{l=0}^d  H_l^0(s-v,v), \quad s \ge v,
\en 
for $H_l^0(\cdot,v)$, $0\le l\le d$, corresponding to the (flattened) branching process~$\hX_0$ of 
Lemma~\ref{AB-dominating-processes},
taken to have initial condition $H_l^0(0,v) = \sum_{j\in \bJ_v} (\l (v-\bt_j))^l/l! \in \s(\tH(v))$,
$0\le l\le d$. Note that $W(v,v) = W^*(v)$.
The error involved in replacing $W^*(s)$ by $W(s,v)$ is bounded in the following lemma.

\begin{lemma}\label{AB-W-to-Wsv}
For $v \le s \le \l^{-1}\log\L$, we have
\[
    \ex\Blb\bigl|\ell(\log[\hc_d e^u W^*(s)] + u) - \ell(\log[\hc_d e^u W(s,v)] + u)\bigr|\, 
        \Giv \tff_v\Brb  \Le    4 \hc_d e^u W^*(v) \etaL \log\L.
\]
\end{lemma}

\proof
We once more use Lemma~\ref{AB-dominating-processes} to justify that both $W^*(s)$ and~$W(s,v)$
belong to the interval
\eq\label{AB-W*-Wsv-sandwich}
    \Bigl[e^{-\l s}\sum_{l=0}^d H_l^-(s-v,v) , e^{-\l s}\sum_{l=0}^d H_l^+(s-v,v)\Bigr],            
\en
where the processes $\hX_-(\cdot,v)$ and~$\hX_+(\cdot,v)$ both have the same initial
condition as~$\hX_0(\cdot,v)$.  Now
$$
    \ex\Blb \sum_{l=0}^d H_l^+(s-v,v) \Giv \tff_v \Brb \Eq e^{\l_+(s-v)} \sum_{l=0}^d H_l^0(0,v)
$$
and
$$
    \ex\Blb \sum_{l=0}^d H_l^-(s-v,v) \Giv \tff_v \Brb \Eq e^{\l_-(s-v)} \sum_{l=0}^d H_l^0(0,v);
$$
hence
\eqs
   \lefteqn{ \ex\bigl\{|W^*(s) - W(s,v)| \giv \tff_v\bigr\} }\\
     &&\Le e^{-\l s}\ex\Blb \sum_{l=0}^d \{H_l^+(s-v,v) -  H_l^-(s-v,v)\} \Giv \tff_v \Brb \\
     &&\Le  e^{-\l v}\sum_{l=0}^d H_l^0(0,v)\{e^{(\l_+ - \l)(s-v)} - e^{(\l_- - \l)(s-v)}\} 
        \Le 4W^*(v) \etaL \log\L,
\ens
by~\Ref{Lambda-conds}.  This, together with \Ref{ADB-l-def} and Lemma~\ref{ADB-phi-smooth}, implies that
\eqs
    \lefteqn{\ex\Blb\bigl|\ell(\log[\hc_d e^u W^*(s)] + u) - \ell(\log[\hc_d e^u W(s,v)] + u)\bigr|\, \Giv \tff_v\Brb }\\
        && \Le   \ex\{ \hc_d e^u |W(s,v) - W^*(s)|\,\giv \tff_v\} \Le 4 \hc_d e^u W^*(v) \etaL \log\L,
\ens
as required.
\ep

We now combine the results of Lemmas \ref{AB-variance-lemma}--\ref{AB-W-to-Wsv} to give the following result,
relating the distribution of $L_{t_\L(u)}/L$ to that of $\ell(\log[\hc_d e^u W(s,v)] + u)$.

\begin{corollary}\label{summary-1}
Take $v := \a_1 \l^{-1}\log\L$ and $s := \a_2 \l^{-1}\log\L$ for $0 < \a_1 < \a_2 < 1$, and fix $0 < \h < \rrr(d)$. 
 Then there is an event $B_{\ref{summary-1}}(\gda,\h,v) \in \tff_v$, and constants $C^0_{\ref{summary-1}}
:= C^0_{\ref{summary-1}}(u_0,d)$, 
$C^1_{\ref{summary-1}}:= C^1_{\ref{summary-1}}(u_0,d)$ and~$C^2_{\ref{summary-1}}:= C^2_{\ref{summary-1}}(d)$, 
such that
\[
   \ex\Blb |f(L_{t_\L(u)}/L) - f(\ell(\log[\hc_d W(s,v)] + u))|  \Giv \tff_v\Brb
          \Le  C^0_{\ref{summary-1}} \|f\|_\infty p_\L + C^1_{\ref{summary-1}} \|f'\|_\infty \e_\L ,
\]
and such that $\pr[\{B_{\ref{summary-1}}(\gda,\h,v)\}^c] \Le C^2_{\ref{summary-1}} p_\L$,
where
\eqs
   \e_\L &:=& \L^{\gda}\{\L^{-\a_2/2}(\log\L)^{d/2} + \L^{\a_2-1}(\log\L)^d + \L^{-\a_1}
                               + \etaL\log\L \} + \L^{-\a_2\h}; \\
    p_\L &:=& \L^{-\gda/2} + \L^{-(\rrr(d)-\h)}\log\L.
\ens
\end{corollary}

\proof
We take the results of Lemmas \ref{AB-variance-lemma}--\ref{AB-W-to-Wsv} in turn.
Using Lemma~\ref{AB-dominating-processes}, we have
\eqa
  \ex\{\l^d \bM_s + \bN_s \giv \tff_v\} &\le& \ex\{ d!H_d^+(s) + H_0^+(s) \giv \tH(v)\} \non\\
    &\le& d! \ex\{W^+(s) e^{\l_+ s} \giv \tH(v)\} \Le 2\, d! W^*(v) e^{\l s}. \label{summ-1}
\ena
Define the event $B_{\ref{summary-1}}\ui(\g,v) := \{W^*(v) \le \L^\gda\}$, 
whose probability is at most $\L^{-\gda}$, by Markov's inequality.
Then, from Lemma~\ref{AB-variance-lemma} and~\Ref{summ-1}, it follows that
\[
 \ex\{\var\{L_t/L \giv \tff_s\} \giv \tff_v\} \Le   C_{\ref{AB-variance-lemma}}  
                  2\, d! W^*(v)\L^{-2}(1 + (\l s)^d) e^{\l(2t-s)}, 
\]
implying that, on $B_{\ref{summary-1}}\ui(\g,v)$,  we have
\eq\label{cor-split-1}
   \ex\Blb |1 - (L_{t_\L(u)}/L) - \ex\{1 - (L_{t_\L(u)}/L) \giv \tff_s\}| \Giv \tff_v\Brb
           \Le C_a(d) \L^{\gda - \a_2/2} e^u (\log\L)^{d/2}.
\en
Next, from Lemma~\ref{AB-first-representation} and~\Ref{summ-1} and on the event $B_{\ref{summary-1}}\ui(\g,v)$, we have
\eqa
  \lefteqn{\ex\Blb\bigl|\ex\{1 - (L_{t_\L(u)}/L) \giv \tff_s\} 
             - \ex\bigl\{\exp\{-\mmm^K_{s,t_\L(u)}\} \giv \tff_s \bigr\}\bigr|  \Giv \tff_v\Brb } \non\\
    &&\Le  C_b(d) W^*(v)(\log\L)^d  \{\L^{-1} e^{\l s}  +  e^{2u} e^{-\l s)} \}  \non\\
   &&\Le C_b(d) \L^\gda (\log\L)^d  \{\L^{\a_2-1}   +  e^{2u} \L^{-\a_1} \}. \label{cor-split-2}
\ena
Turning to Lemma~\ref{AB-M-approx}, we find that
\eqa
     \lefteqn{ \ex\Blb \bigl|\ex\bigl\{\exp\{-\mmm^K_{s,t_\L(u)}\} \giv \tff_s \bigr\} 
           - \ex\bigl\{\exp\{-\mmc^K_{s,t_\L(u)}\} \giv \tff_s \bigr\}\bigr|\,
                                     I[B_{\ref{AB-M-approx}}(\gda,s)] \Giv \tff_v \Brb }\non\\
    &&\Le C_c(d) \L^{\gda} e^u \{\L^{-1} \{\log\L\}^d e^{\l s} + \etaL\}\phantom{XXXXXXXXXXX}\non \\
    &&\Eq C_c(d) \L^{\gda} e^u \{\L^{\a_2-1} \{\log\L\}^d  + \etaL\}. \label{cor-split-3}
\ena
Then, from Lemma~\ref{AB-e-M-approx}, we have
\eqa
 \lefteqn{ \ex\Blb \bigl|\ex\{e^{-\mmc^K_{s,t_\L(u)}} \giv \tff_s\}  - \ell(\log[\hc_d W^*(s)] + u)\bigr|\,
         I[B_{\ref{AB-e-M-approx}}(\gda,\h,s)] \Giv \tff_v\Brb }\non\\
     &&\Le C'_{\ref{AB-e-M-approx}}(1+e^u)(\L^{\gda}(\etaL\log \L + \L^{-1}e^{\l s}) + e^{-\l\h s}) \non\\
     &&\Eq C'_{\ref{AB-e-M-approx}}(1+e^u)(\L^{\gda}(\etaL\log \L + \L^{\a_2-1}) + \L^{-\a_2\h}). \label{cor-split-4}
\ena
Finally, from Lemma~\ref{AB-W-to-Wsv}, on the event $B_{\ref{summary-1}}\ui(\g,v)$, we have
\eq\label{cor-split-5}
    \ex\Blb\bigl|\ell(\log[\hc_d W^*(s)] + u) - \ell(\log[\hc_d W(s,v)] + u)\bigr| \Giv \tff_v\Brb
         \Le    4 \hc_d e^u \L^{\gda} \etaL \log\L.
\en
Combining \Ref{cor-split-1} to~\Ref{cor-split-5}, we deduce that, on the event $B_{\ref{summary-1}}\ui(\g,v)$,
and uniformly in $u \le u_0$,
\eqa
    \lefteqn{ \ex\Blb \bigl|(L_{t_\L(u)}/L) - \ell(\log[\hc_d W(s,v)] + u)\bigr|\, 
            I[\hB(\gda,\h,s)]\Giv \tff_v\Brb }.  \non\\
       &&\Le C_*(d,u_0) \bigl(\L^{\gda}\{\L^{-\a_2/2}(\log\L)^{d/2} + \L^{\a_2-1}(\log\L)^d + \L^{-\a_1}
                               + \etaL\log\L \} + \L^{-\a_2\h}\bigr) \non\\
       &&\ =: C_*(d,u_0)\e_\L,            \label{cor-split-*}
\ena
where $\hB(\gda,\h,s) := B_{\ref{AB-e-M-approx}}(\gda,\h,s) \cap B_{\ref{AB-M-approx}}(\gda,s)$.

For the exceptional set, from Lemmas \ref{AB-e-M-approx} and~\ref{AB-M-approx}, we have
\eqs
     \pr[\hB(\gda,\h,s)\}^c]
     &\le&  C_{\ref{AB-e-M-approx}}\{\L^{-\gda} + \l s e^{-2\l (\rrr(d)-\h)s}) + C_{\ref{AB-M-approx}}\L^{-\gda}\} \\
     &\le& C_e(d) \{\L^{-\gda} + \L^{-2(\rrr(d)-\h)}\log\L\}.
\ens
On the other hand, for any set~$B \in \ff$ with $\pr[B]=p$, and for any $\s$-field $\cG \subset \ff$, we have
\[
      p \Eq \pr[B] \ \ge\ \pr[\{\pr[B\giv\cG] > \sqrt p\}] \sqrt p,
\]
by the total probability formula, implying that $\pr[B\giv\cG] \le \sqrt p$ with probability at least $1-\sqrt p$.
Hence there is an event $B_{\ref{summary-1}}\ut(\gda,\h,v) \in \tff_v$, whose complement has probability at most
\eq\label{excep-prob}
     (C_e(d))^{1/2} \{\L^{-\gda/2} + \L^{-(\rrr(d)-\h)}\log\L\}\ =:\ (C_e(d))^{1/2} p_\L,
\en
on which $\pr[\{\hB(\gda,\h,s)\}^c \giv \tff_v] \le (C_e(d))^{1/2} p_\L$.  Now define
$Z_u := \ell(\log[\hc_d W(s,v)] + u)$ and $Y_u := L_{t_\L(u)}/L$.  Then, for any bounded Lipschitz function~$f$,
we conclude from \Ref{cor-split-*} and~\Ref{excep-prob} that, for $u \le u_0$ and on the event
\[
    B_{\ref{summary-1}}(\gda,\h,v) \Def B_{\ref{summary-1}}\ui(\gda,v) \cap B_{\ref{summary-1}}\ut(\gda,\h,v),
\]
we have
\eqs
   \lefteqn{ \ex\Blb |\ex\{f(Y_u)\} - \ex\{f(Z_u)\}|\,\giv \tff_v \Brb }\\
     &&\Le  \ex\Blb |\ex\{f(Y_u)\} - \ex\{f(Z_u)\}|I[\hB(\gda,\h,s)] \right.\\
      &&\qquad\qquad\mbox{} \left.      +  |\ex\{f(Y_u)\} - \ex\{f(Z_u)\}|I[\{\hB(\gda,\h,s)\}^c]\,\Giv \tff_v \Brb \\
     &&\Le \|f'\|_\infty \ex\{|Y_u - Z_u|I[\hB(\gda,\h,s)] \giv \tff_v\} 
             + 2\|f\|_\infty \pr[\{\hB(\gda,\h,s)\}^c \giv \tff_v] \\
     &&\Le \|f'\|_\infty C_*(d,u_0) \e_\L + 2\|f\|_\infty (C_e(d))^{1/2} p_\L.
\ens
This proves the corollary.
\ep

\subsection{The main theorem}
We now use Corollary~\ref{summary-1} to compare the conditional distributions, given~$\tff_v$, of the
{\it normalized\/} random variables $Y(u,v)$ and~$Z(u,v)$, 
where
\eqa
  Y(u,v) &:=& e^{\l v/2}\{(L_{t_\L(u)}/L) - \ell(\log[\hc_d W^*(v)] + u)\}; \non\\
  Z(u,v) &:=& e^{\l v/2}\{\ell(\log[\hc_d W(s,v)] + u) - \ell(\log[\hc_d W^*(v)] + u)\}, \label{YandZ-def}
\ena
for a careful choice of~$s$, with the centring constant $\ell(\log[\hc_d W^*(v)] + u)$
chosen because $W^*(v) = \ex\{W(s,v) \giv \tff_v\}$.  These are the correct standardizations to achieve a non-trivial 
limit.  Thus we wish to compare $\ex f(Y(u,v))$ with $\ex f(Z(u,v))$, for Lipschitz functions~$f$ that 
have $\|f\|_\infty \le 1$ and $\|f'\|_\infty \le 1$.  This corresponds to taking $\|f\|_\infty \le 1$
and $\|f'\|_\infty \le e^{\l v/2}$ in Corollary~\ref{summary-1}, because of the pre-factors $e^{\l v/2}$
in the definitions of $Y(u,v)$ and~$Z(u,v)$.  Thus, although $p_\L$ is already small for large~$\L$, 
if $\h < \rrr(d)$, we need
also to show that, for $v = \a_1\l^{-1}\log\L$, it is possible to choose $\a_2$, $\h$ and~$\g$ so as to
make $e^{\l v/2}\e_\L = \L^{\a_1/2}\e_\L$ small with~$\L$. Recalling the definition~\Ref{AB-eta-Lambda-def} 
of~$\etaL$, the expression for~$\e_\L$ in Corollary~\ref{summary-1} shows that this
is the case, for~$\gda>0$ chosen small enough, if,  
$$
       \a_1 \ <\ \a_2;\quad  \a_2 \ <\ 1 - \a_1/2 ;\quad \gda < \a_1/2;\quad \a_1 < 2\a_2\h
     \ \mbox{ and }  \a_1\ <\ 2\g_g/d .
$$
So, for
\[
   \a_1 \ <\ 2\min\{\g_g/d,\rrr(d)/(1 + \rrr(d)\},
\]
choose $0 < \h < \rrr(d)$ so that $2\h/(1+\h) > \a_1$ and then~$\a_2$ so that $\a_1/(2\h) < \a_2 < 1-\a_1/2$;
then, if we choose
\eq\label{ADB-gda-choice}
   0\ <\ \gda \Eq \tfrac23\min\{\g_g/d - \a_1/2, (\a_2-\a_1)/2,1-\a_1/2-\a_2,\a_1/2,\a_2\h-\a_1/2\},
\en
it follows that there are constants~$C = C(d,u_0)$ and $C' = C'(d)$ such that
\eq\label{Y-Z-bound}
    |\ex\{f(Y(u,v)) \giv \ff_v\} - \ex\{f(Z(u,v)) \giv \ff_v\}| 
           \Le C\{\L^{-\gda/2}(\log\L)^d  + \L^{-(\rrr(d)-\h)}\},
\en
for all $f \in \Fbw$, except on an event of probability at most $C'\{\L^{-\gda/2}  + \L^{-(\rrr(d)-\h)}\}$.
Particular choices are to take 
\eq\label{eta-and-alpha2-choice}
     \h \Def \frac12\Bl \rrr(d) + \frac{\a_1}{2-\a_1} \Br,\quad\mbox{and}\quad 
    \a_2 \Def \frac12\Blb 1 + \frac{\a_1}2 \Bl \frac1{\h}-1 \Br \Brb,
\en
in which case we can take any~$0 < \gda' < \min\{\gda/2,(\rrr(d)-\h)\}$, and express the error
in~\Ref{Y-Z-bound} as $C\L^{-\gda'}$, except on an event of probability at most $C'\L^{-\gda'}$,
albeit with different constants $C=C(u_0,d)$ and~$C'(d)$.

Corollary~\ref{summary-1} and~\Ref{Y-Z-bound} compare the distribution of $L_{t_\L(u)}/L$ with 
that of the quantity $\ell(\log[\hc_d W(s,v)] + u)$,
for any~$u \le u_0$.  The path of~$L_{t_\L(u)}/L$ is approximated, to first order, by a time shift of
the deterministic path~$\ell(u)$, and the shift is the same throughout the path, being determined by
the value of the single $\tff_s$-measurable random variable~$W(s,v)$.  In the remaining argument, we exploit this
to show that, to a good approximation, the path after time~$v$ is that of the approximation 
$\ell(\log[\hc_d W^*(v)] + \cdot)$, together with a perturbation that can be expressed in the form 
$e^{-\l v/2}N h_v(\cdot)$, where $h_v(\cdot)$ is an $\tff_v$-measurable function depending on the value of~$W^*(v)$,
and~$\law(N \giv \tff_v)$ is the standard normal distribution.

To do so, in view of~\Ref{Y-Z-bound}, we now need a central limit theorem for $Z(u,v)$ as defined in~\Ref{YandZ-def}.
Writing 
\eqa
   K_{2}(u,v) &:=& (D\ell)(u + \log[\hc_d W^*(v)])/W^*(v) \Eq k \frac d{dx}\{\ell(\log x)\}\Bigr|_{k W^*(v)},\phantom{XX}
             \label{ADB-Kuv-defs}
\ena
where the final equality holds for all $k > 0$, the next lemma shows that $Z(u,v)$ is close in distribution to
$K_2(u,v)\,e^{\l v/2}\{W(s,v) - W^*(v)\}$.

\begin{lemma}\label{K-version}
Let $Z(u,v)$ be defined as in~\Ref{YandZ-def}, and let $v := \a_1\l^{-1}\log\L$ and
$s := \a_2\l^{-1}\log\L$;  suppose that $\gda$ is as for \Ref{ADB-gda-choice} and~$\gda' = \half\min\{\gda/2,(\rrr(d)-\h)\}$,
where~$\h$ is as in~\Ref{eta-and-alpha2-choice}.  
Then there is a constant~$C = C(d,u_0)$ such that, for all $f \in \Fbw$, and on the event $\{W^*(v) \le \L^\gda\}$,
\[
     |\ex\{f(Z(u,v)) \giv \tff_v\} - \ex\{f(K_2(u,v)\,e^{\l v/2}\{W(s,v) - W^*(v)\})\} \giv \tff_v\}|
              \Le C\L^{-\gda'},
\]
uniformly in $u \le u_0$. 
\end{lemma}

\proof
From~\Ref{ADB-l-def}, we have $g(x) := \ell(\log x) = 1 - \ex\{e^{-xW}\}$, so that,
by Taylor's expansion, for any $x,y > 0$, we can write
\[
    |g(x+y) - (g(x) + yg'(x))| \Le \half y^2 \|g''\|_\infty \Eq \half y^2 \ex W^2 \Le \half y^2.
\]
from~\Ref{AB-MG-variance}. Thus, in making a linear approximation to
$$
    \ell(\log[k W(s,v)]) - \ell(\log[k W^*(v)]) \Eq g(k W(s,v)) - g(kW^*(v)), 
$$
the remainder term can be bounded by $\half k^2 (W(s,v)-W^*(v))^2$.  Now, because $W^*(v) = \ex\{W(s,v) \giv \tff_v\}$,
we have
$$
    \ex\{ (W(s,v)-W^*(v))^2 \giv \tff_v\} \Eq V(s,v) \Def \var(W(s,v) \giv \tff_v) \Le W^*(v) e^{-\l v},
$$ 
where the inequality follows using~\Ref{AB-MG-variance}.
Hence, for any $k > 0$, and using~\Ref{ADB-Kuv-defs}, we have
\eqa
     \lefteqn{\ex\bigl\{\bigl| e^{\l v/2}\{\ell(\log[k W(s,v)]) - \ell(\log[k W^*(v)])\}
        - e^{\l v/2} \{W(s,v) - W^*(v)\}\,K_2(u,v)\bigr| \giv \tff_v\bigr\} } \non\\
              && \Le  \half k^2 V(s,v) \Le \half k^2 W^*(v) e^{-\l v}.  
                  \phantom{XXXXXXXXXXXXXXXXXXXXXXX} \label{Taylor-0}
\ena
Thus, taking $k = \hc_d e^{u}$ in~\Ref{Taylor-0}, and on~$\{W^*(v) \le \L^\gda\}$, it follows that
\eqa
  \lefteqn{ 
    \ex\Bigl\{\bigl| e^{\l v/2} \{\ell(\log[\hc_d W(s,v)]+u) - \ell(\log[\hc_d W^*(v)]+u)\} } \non \\
    &&\qquad\qquad\mbox{}  -
         K_{2}(u,v) e^{\l v/2}\{W(s,v) - W^*(v)\} \bigr| \,\Big|\, \tff_v\Bigr\}  \non \\ 
     &&\qquad\Le \half \hc_d^2 e^{2u} \L^{\gda - \a_1/2}, 
             \phantom{XXXXXXXXXXXXXXXXXXXX} 
      \label{AB-Taylor}
\ena
and the lemma follows because $\gda' + \gda < 3\gda/2 \le \a_1/2$, from~\Ref{ADB-gda-choice}.
\ep

We are now in a position to prove a central limit theorem, with an error bound expressed in terms of the
bounded Wasserstein distance.

\begin{theorem}\label{AB-main-theorem}
Suppose that $v = \a\l^{-1}\log\L$ for $0 < \a < 2\min\{\g_g/d,\rrr(d)/(1 + \rrr(d)\}$,
where~$\g_g$ is as in~\Ref{AB-volume-approx} and~$\rrr(d)$ as in~\Ref{AB-min-r-eta-def}
(so that $\rrr(d) = 1/2$ for $d \le 6$).  Suppose that $\gda$ is as for \Ref{ADB-gda-choice} 
and~$\gda' = \half\min\{\gda/2,(\rrr(d)-\h'),(\a_2-\a)\}$,
where~$\h'$ and~$\a_2$ are as in~\Ref{eta-and-alpha2-choice}, with $\a_1=\a$.
Suppose that~$\L$ is large enough that~\Ref{Lambda-conds} is satisfied, and that~$\L^{4\a\rrr(d)/7} \ge (d+1)^3$
and $\a\log\L > c_{1*}$, where~$c_{1*}$ is as in Theorem~\ref{AB-branching-appxn}.
Then, for any $u_1 < u_0 \in \re$, there exist constants $C(d,u_1,u_0)$ and~$C'(d,u_1,u_0)$ 
and an event $E^*(v) \in \s(\tH(v))$ with $\pr[E^*(v)^c] \le C'(d,u_1,u_0)\L^{-\gda'}$ such that
\eqs
    \lefteqn{\dbw\bigl(\law\{e^{\l v/2}((L_{t_\L(u)}/L) - \ell(\log[\hc_d W^*(v)] + u)) \giv \tff_v \cap E^*(v)\},
             \nn(0,\{K_2(u,v)\}^2 W^*(v) / (d+1)) \bigr) } \\
       && \Le C(d,u_1,u_0)\L^{-\gda'}, \phantom{XXXXXXXXXXXXXXXXXXXXXXXXXXXXXXXXX}
\ens
uniformly in $u_1 \le u \le u_0$, where~$K_2(u,v)$ is defined in~\Ref{ADB-Kuv-defs}, $\hc_d$ in Lemma~\ref{AB-e-M-approx}
and $t_\L(u)$ in~\Ref{AB-t-Lambda-def}.
\end{theorem}

\proof
In view of~\Ref{Y-Z-bound} and Lemma~\ref{K-version}, it suffices to show that
\[
   \dbw\bigl(\law(e^{\l v/2}\{W(s,v) - W^*(v)\} \giv \tff_v),\nn(0, W^*(v) / (d+1))\bigr) \Le C_1(d,u_1,u_0)\L^{-\gda'},
\] 
with $s = \a_2\l^{-1}\log\L$ and~$\a_2$ as in~\Ref{eta-and-alpha2-choice}.
Corollary~\ref{AB-CLT-cor}, with $\h = 6\rrr(d)/7$,  shows that there is an event~$E^\h(v) \in \tH(v)$ with
$\pr[\{E^\h(v)\}^c] \le C'(d)\L^{-2\a\rrr(d)/7}$ such that, on~$E^\h(v)$, 
\eqs
   \lefteqn{ \dbw\bigl( \law(e^{\l v/2}\{W(s,v) - W^*(v)\} \giv \tff_v),
                  \law(e^{\l v/2}X_v\uh(t_v^{-1}(s-v) \giv \tff_v)\bigr) }\\
           &&\Le C(d)\{L^{-\a\rrr(d)/7}\}, \phantom{XXXXXXXXXXXXXXX}
\ens
provided that $\L^{4\a\rrr(d)/7} \ge (d+1)^3$.
Then, from \Ref{AB-t(s)-def} and~\Ref{AB-brownian-repn}, 
\[
   \law\bigl(e^{\l v/2} X_v\uh(t_v^{-1}(s-v))\bigr) \Eq \nn\Bigl(0,\frac{W^*(v)}{d+1}(1-e^{-\l(s-v)})\Bigr),
\]
and the theorem follows because $\dbw(\nn(0,\s_1^2),\nn(0,\s_2^2)) = O(|\s_1 - \s_2|)$ 
and 
$$
         W^*(v)e^{-\l(s-v)} \Eq W^*(v)\L^{-(\a_2-\a)} \Le L^{\gda' -(\a_2-\a)},
$$
on $\{W^*(v) \le \L^{\gda'}\}$, and $\gda' \le \half(\a_2-\a)$, from~\Ref{ADB-gda-choice}.
\ep

This theorem is not quite the same as Theorem~\ref{ADB-main-theorem}, because both mean and
variance are expressed in terms of~$W^*(v) = W(v,v)$, which, as is seen from its definition in~\Ref{ADB-W-star-def},
is not necessarily determined by knowledge
of~$\law_v$ alone, because all the birth times of~$\bX$ come into its definition.  Instead,
one can observe~$\hW(v)$ as in~\Ref{ADB-W-star-def-new}.  We now show that this is enough. 

We construct a lower bound $\hW_-(v)$ for~$\hW(v)$ by summing over the subset of the birth times~$\hJ_v \subset J_v$
in~\Ref{ADB-W-star-def-new} that belong to $J_v \cap \tJ_v$, where~$J_v$ is defined in~\Ref{ADB-Js-etc-defs}, and
\[
    \tJ_v \Def \Bigl\{ j\ge0\colon \bP_j \notin \bigcup_{{l\in\bJ_v}\atop{l < j}} \KK(\bP_j,2v)\Bigr\},
\]
with $\bJ_v$ the birth times of~$\bX$ before~$v$, defined in~\Ref{ADB-bJ-etc-defs}.
These give rise to non-intersecting neighbourhoods at time~$v$, though not 
necessarily to all such, and they form a subset more amenable to calculation.  Then
it is immediate from~\Ref{AB-volume-approx} that, for all~$\L$ sufficiently large,
\eqs
   \ex|\bJ_v \setminus \tJ_v| &\le& 2\ex\{|\bJ_v|(|\bJ_v|-1)L^{-1}(2v)^d\vva \} \\
              &\le& 2c_2 e^{2\l_+ v}\L^{-1}(2\l v)^d,
\ens
the final inequality following from~\Ref{AB-N-star-moments}.  Then, using arguments analogous to those
in Lemma~\ref{AB-variance-lemma}, we have
\eqs
    \ex|\bJ_v \setminus J_v| &\le& \ex\Bigl\{ \int_0^v L^{-1}\vva  M_d^+(u)\,c_1 e^{\l_+(v-u)} M_0^+(du) \Bigr\} \\
       &=& \ex\Bigl\{ \r\vva \int_0^v L^{-1}\vva  (M_d^+(u))^2\,c_1 e^{\l_+(v-u)} \,du \Bigr\} \\
       &\le& \r\vva ^2 L^{-1} c_1(c_2 d! \l_+^{-d})^2 \int_0^v e^{\l_+(v+u)}\,du \Le C\L^{-1}e^{2\l_+ v}.
\ens  
Hence, for $v \le \tlm$,
\[
    0 \Le \ex\{\hW(v) - \hW_-(v)\} \Eq O\Blb \L^{-1}\sum_{l=0}^d e^{\l v}(\log\L)^{d+l} \Brb,
\]
and, for~$v = \a\l^{-1}\log\L$, this is of order $O(\L^{-1+\a}(\log\L)^{2d})$.  The most
sensitive place where this enters is into~$\ell(\log[\hc_d W^*(v)] + u)$, when the difference has to be small
relative to $\L^{-\a/2}$, because of the factor~$e^{\l v/2}$; but this is the case
if $\a < 2/3$, as in the statement of the theorem, by Lemma~\ref{ADB-phi-smooth}.  The conversion of~$E^*(v)$ into an
event that can be determined from~$\law_v$ can be accomplished in similar fashion, by
modifying the definitions of its constituent events in terms of $W_j(v)$, $0\le j\le v$.

\section*{Appendix}
\setcounter{equation}0
\setcounter{section}4
\setcounter{theorem}0
We note here two technical lemmas that are used in the previous arguments.
The first establishes a bound on the extreme fluctuations of an integral with respect to a
compensated Poisson process.

\begin{lemma}\label{AB-process-bounds}
  Let $X(t) := \int_0^t F(u)\{Z(du) - du\}$, where~$Z$ is a Poisson process and the process~$F$
is predictable and a.s.\ bounded in modulus by the deterministic function~$G$.  Define
$G_2(s,t) := \int_s^t \{G(u)\}^2\,du$ and $G^*(s,t) := \sup_{s\le u\le t}G(u)$. Then
\[
   \pr\Bigl[ \sup_{t_1\le t\le t_2}|X(t) - X(t_1)| > a\Bigr]
    \Le 2\exp\bigl\{-a^2/\{2eG_2(t_1,t_2)\}\bigr\},
\]
for all $0 \le a \le eG_2(t_1,t_2)/G^*(t_1,t_2)$.  If~$G$ is decreasing, we have
\[
   \pr\Bigl[ \sup_{t_1\le t\le t_2}|X(t) - X(t_1)| > a\Bigr]
    \Le 2\exp\bigl\{-a^2/\{2e\{G(t_1)\}^2(t_2-t_1)\}\bigr\},
\]
for all $0 \le a \le eG(t_1)(t_2-t_1)$.
\end{lemma}

\proof
For any $\th$, the process
\[
   Y(t) \Def \exp\Blb \th X(t) - \int_0^t\{e^{\th F(u)} - 1 - \th F(u)\}\,du \Brb
\]
is a supermartingale (van de Geer~(1995, p.~1795)), and stopping at~$a$ easily yields
\eqs
    \pr\Bigl[ \sup_{t_1\le t\le t_2} (X(t) - X(t_1)) > a\Bigr] &\le&
      e^{-\th a} \ex\Blb \exp\Bl \int_{t_1}^{t_2} \{e^{\th F(u)} - 1 - \th F(u)\}\,du \Br \Brb\\
   &\le& e^{-\th a} \exp\Bl \frac e2  \th^2 G_2(t_1,t_2)\,du \Br,
\ens
if $0 \le \th G^*(t_1,t_2) \le 1$.  The corresponding bound for $\inf_{t_1\le t\le t_2} (X(t) - X(t_1))$
is proved in analogous fashion.
Now, if $a \le eG_2(t_1,t_2)$, choose $\th = a/\{eG_2(t_1,t_2)\}$,
giving the first conclusion of the lemma.  The second follows by choosing 
$\th = a/\{e G(t_1)^2 (t_2-t_1)\}$.
\ep

\medskip
The second lemma establishes some smoothness of the function~ $\f^1_s(\th) := \ex\{e^{-\th W^1(s)}\}$.

\begin{lemma}\label{ADB-phi-smooth}
 With~$\f^1_s$ defined as above, and for any $s,h,\th > 0$, we have
\eqs
   &&|\f^1_{s+h}(\th) - \f^1_s(\th)| \Le \half\th e^{-1} e^{-s}(1-e^{-h});\\
   &&|\f^1_{s}(\th(1+\d)) - \f^1_s(\th)| \Le    \d \min\{e^{-1},\th\}.
\ens
\end{lemma}

\proof
We note that $W^1(s) \ge 0$ and that $\ex W^1(s) = 1$ for all~$s$.  Then, writing 
$X_s(h) := W^1(s+h) - W^1(s)$ and using~\Ref{AB-MG-variance}, we have
\eq\label{ADB-easy-facts}
    \ex\{X_s(h) \giv \hff_s\} \Eq 0;\qquad \ex\{(X_s(h))^2 \giv \hff_s\} \Le W^1(s) e^{-s}(1-e^{-h}),
\en
for any $s,h > 0$. Hence, using~\Ref{ADB-easy-facts}, and taking expectations first conditional on~$\hff_s$, we have
\eqs
   \f^1_{s+h}(\th) - \f^1_s(\th) &=& \ex\bigl\{e^{-\th W^1(s)}\{(e^{-\th X_s(h)}-1+\th X_s(h)) - \th X_s(h)\}\bigr\} \\
      &=& \ex\bigl\{e^{-\th W^1(s)}\ex\{(e^{-\th X_s(h)}-1+\th X_s(h)) \giv \hff_s\} \bigr\}.
\ens 
This implies that
\[
   | \f^1_{s+h}(\th) - \f^1_s(\th)| \Le \half\ex\{e^{-\th W^1(s)}\th^2 W^1(s) e^{-s}(1-e^{-h})\}
           \Le \frac\th{2e} e^{-s}(1-e^{-h}),
\]
since $xe^{-x} \le e^{-1}$, proving the first inequality.

For the second, since $e^{-x}(1-e^{-\d x}) \le \d e^{-1}$ in $x\ge0$
and $\ex W^1(s)=1$,
\eqs
  \qquad\quad|\f^1_{s}(\th(1+\d)) - \f^1_s(\th)| &=& |\ex\{e^{-\th W^1(s)}(1 - e^{-\th\d W^1(s)})\}|
                \Le \d \min\{e^{-1}, \th\}. \hskip0.5in\halmos
\ens

\section*{Acknowledgement}
ADB thanks the Department of Statistics and Applied Probability at the National University of Singapore,
and the mathematics departments of the University of Melbourne and Monash University, for
their kind hospitality while much of the work was undertaken. 
AR thanks the School of Mathematics and Statistics at the University of Melbourne for their kind hospitality.

\end{document}